\documentclass{gtart}

\def\ifplaintex{\expandafter\ifx\csname documentclass\endcsname\relax}

\def\gtp{{\mathsurround=0pt\it $\cal G\mskip-2mu$eometry \&\ 
$\cal T\!\!$opology $\cal P\!$ublications}}  

\def\recd{{\small Received:\qua\receiveddate\ifx\reviseddate\relax
\else\qquad Revised:\qua\reviseddate\fi\par}} 

\def\lognumber#1{\def\thelognumber{#1}}
\def\volumenumber#1{\def\thevolumenumber{#1}}
\def\volumeyear#1{\def\thevolumeyear{#1}}
\def\papernumber#1{\def\thepapernumber{#1}}
\def\pagenumbers#1#2{\def\startpage{#1}\def\finishpage{#2}}
\def\published#1{\def\publishdate{#1}}

\def\received#1{\def\receiveddate{#1}}

\def\accepted#1{\def\accepteddate{#1}}

\long\def\asciiabstract#1{\long\def\theasciiabstract{#1}}
\def\asciikeywords#1{\def\theasciikeywords{#1}}

\let\\\par\let\thelognumber\relax\let\thevolumenumber\relax
\let\thepapernumber\relax\let\thevolumeyear\relax\let\startpage\relax
\let\finishpage\relax\let\publishdate\relax\let\receiveddate\relax
\let\reviseddate\relax\let\accepteddate\relax\let\theasciititle\relax
\let\theasciiauthors\relax
\let\theasciiabstract\relax\let\theasciikeywords\relax

\let\theasciiemail\relax

\ifplaintex
\font\logobig=cmssbx10 scaled 3836
\font\logomed=cmssbx10 scaled 2557
\else
\font\logobig=cmssbx10 scaled 4200
\font\logomed=cmssbx10 scaled 2800
\fi

\long\def\makeagttitle{   
\count0=\startpage
\agt\hfill      
\hbox to 45truept{\vbox to 0pt{\vglue -13truept{\logomed A\kern -.37em{\logobig 
T}\kern -.38em G}\vss}\hss}
\break
{\small Volume \thevolumenumber\ (\thevolumeyear)
\startpage--\finishpage\nl
Published: \publishdate}

\vglue .25truein

{\parskip=0pt\leftskip 0pt plus
1fil\def\\{\par\smallskip}{\Large\bf\thetitle}\par\medskip} \vglue
0.05truein

{\parskip=0pt\leftskip 0pt plus 1fil\def\\{\par}{\sc\theauthors}
\par\medskip}%
 
\vglue 0.03truein 

{\small\leftskip 25truept\rightskip 25truept{\bf Abstract}\stdspace\theabstract

{\bf AMS Classification}\stdspace\theprimaryclass
\ifx\thesecondaryclass\relax\else; \thesecondaryclass\fi\par
{\bf Keywords}\stdspace \thekeywords\par}\vglue 7truept

}   

\ifplaintex
\font\phead=cmsl9 scaled 950
\font\pnum=cmbx10 scaled 913
\font\pfoot=cmsl9 scaled 950
\headline{\vbox to 0pt{\vskip -4.5mm\line{\small\phead\ifnum
\count0=\startpage ISSN 1472-2739 (on-line) 1472-2747 (printed)
\hfill {\pnum\folio}\else\ifodd\count0\def\\{ }%
\ifx\theshorttitle\relax\thetitle\else\theshorttitle\fi\hfill{\pnum\folio}
\else\def\\{ and }{\pnum\folio}\hfill\ifx\theshortauthors\relax\theauthors
\else\theshortauthors\fi\fi\fi}\vss}}
\footline{\vbox to 0pt{\vglue 0mm\line{\small\pfoot\ifnum\count0=\startpage
\copyright\ \gtp\hfill\else
\agt, Volume \thevolumenumber\ (\thevolumeyear)\hfill\fi}\vss}}
\else
\font=cmsl9 scaled 1050
\font\lnum=cmbx10 
\font=cmsl9 scaled 1050
\makeatletter
\def\@oddhead{{\small\lhead\ifnum\count0=\startpage ISSN 1472-2739 
(on-line) 1472-2747 (printed)\hfill {\lnum\number\count0}\else\ifodd\count0
\def\\{ }\ifx\theshorttitle\relax \thetitle \else\theshorttitle\fi\hfill
{\lnum\number\count0}\else\def\\{ and }{\lnum\number\count0}
\hfill\ifx\theshortauthors\relax 
\theauthors\else\theshortauthors\fi\fi\fi}}\def\@evenhead{\@oddhead}
\def\@oddfoot{\small\lfoot\ifnum\count0=\startpage\copyright\ \gtp\hfill\else
\agt, Volume \thevolumenumber\ (\thevolumeyear)\hfill\fi}
\def\@evenfoot{\@oddfoot}
\makeatother
\fi
\let\maketitlepage\makeagttitle
\let\makeshorttitle\maketitlepage
\let\maketitle\maketitlepage

\newwrite\gtoutfile
\long\gdef\makeheadfile{  
{\def\\{, }\def\s{ }
\immediate\openout\gtoutfile head.xxx
\immediate\write\gtoutfile{To: math@arxiv.org}
\immediate\write\gtoutfile{Subject: put OR rep NNNNN:ppppp}
\immediate\write\gtoutfile{--text follows this line--}
\immediate\write\gtoutfile{Proxy-for: \ifx\theasciiauthors\relax
\theauthors\else\theasciiauthors\fi\s<\ifx\theasciiemail\relax\theemail\else\theasciiemail\fi>}
\immediate\write\gtoutfile{\noexpand\\}
\immediate\write\gtoutfile{Authors: \ifx\theasciiauthors\relax
\theauthors\else\theasciiauthors\fi}
{\def\\{ }\immediate\write\gtoutfile{Title: \ifx\theasciititle\relax
\thetitle\else\theasciititle\fi}}
\immediate\write\gtoutfile{Subj-class: GT or SG, GR etc}
\immediate\write\gtoutfile{MSC-class: \theprimaryclass\ifx\thesecondaryclass\relax\else, \thesecondaryclass\fi}
\immediate\write\gtoutfile{Journal-ref: Algebr. Geom. Topol. \thevolumenumber\s
(\thevolumeyear) \startpage-\finishpage}
\immediate\write\gtoutfile{Comments: Published by Algebraic and
Geometric Topology at}
\immediate\write\gtoutfile{\s\s\s  http://www.maths.warwick.ac.uk/agt/AGTVol\thevolumenumber/agt-\thevolumenumber-\thepapernumber.abs.html}
\immediate\write\gtoutfile{\noexpand\\}
\immediate\write\gtoutfile{}
\ifx\theasciiabstract\relax
\immediate\write\gtoutfile{\theabstract}\else
\immediate\write\gtoutfile{\theasciiabstract}\fi
\immediate\write\gtoutfile{}
\immediate\write\gtoutfile{\noexpand\\}
\immediate\write\gtoutfile{}
\immediate\closeout\gtoutfile}}  

\def\maketitlepage{\makeagttitle\makeheadfile}
\let\makeshorttitle\maketitlepage
\let\maketitle\maketitlepage

\lognumber{34}
\volumenumber{3}
\volumeyear{2003}
\papernumber{34}
\published{11 October 2003}
\pagenumbers{1005}{1035}
\received{20 Febuary 2003}
\accepted{8 October 2003}

\usepackage{amsmath,amssymb}
\usepackage{latexsym,amscd}
\usepackage{epsfig}
\theoremstyle{plain}
\newtheorem{theorem}{Theorem}

\newtheorem{lemma}{Lemma}
\newtheorem*{ztor}{$\mathbb{Z}$-torsion theorem}
\theoremstyle{definition}
\newtheorem{definition}{Definition}
\newtheorem{example}{Example}
\newtheorem{remark}{Remark}
\newtheorem{conjecture}{Conjecture}
\newenvironment{Relax}{\relax}{\relax}

\begin{document}

\title{Deformation of string topology into\\homotopy skein modules}
\author{Uwe Kaiser}
\address{Department of Mathematics,
Boise State University\\
1910 University Drive,
Boise, ID 83725-1555, USA}
\email{kaiser@math.boisestate.edu}
\url{http://diamond.boisestate.edu/\char'176kaiser/}
\begin{abstract} 
Relations between the string topology of Chas and Sullivan and the homotopy skein modules of Hoste and Przytycki are studied. This provides new insight into the structure of homotopy skein modules and their meaning in the framework of quantum topology. Our results can be considered as \textit{weak extensions} to all orientable $3$-manifolds of classical results by Turaev and Goldman concerning intersection and skein theory on oriented surfaces.
\end{abstract}
\asciiabstract{Relations between the string topology of Chas and
Sullivan and the homotopy skein modules of Hoste and Przytycki are
studied. This provides new insight into the structure of homotopy
skein modules and their meaning in the framework of quantum
topology. Our results can be considered as weak extensions to all
orientable 3-manifolds of classical results by Turaev and Goldman
concerning intersection and skein theory on oriented surfaces.}
\begin{Relax}\end{Relax}

\primaryclass{57M25}
\secondaryclass{57M35, 57R42}
\keywords{$3$-manifold, string topology, deformation, skein module, torsion, link homotopy, free loop space, Lie algebra}
\asciikeywords{3-manifold, string topology, deformation, skein module, torsion, link homotopy, free loop space, Lie algebra}
\maketitle 

\setcounter{section}{-1}

\section{Introduction}

In 1999 Moira Chas and Dennis Sullivan discovered the structure of a
graded Lie algebra on the equivariant homology of the free loop space
of an oriented smooth or combinatorial $d$-manifold \cite{CS1}.
Later Cattaneo, Fr\"{o}hlich and Pedrini developed ideas about the
quantization of string topology
in the framework of topological field theory \cite{CFP}.
It is the goal of this paper to study Vassiliev-Kontsevitch 
and skein theory of links in $3$-manifolds
in relation with string topology. 
 Our approach is \textit{intrinsically} $3$-dimensional. It hints towards a general deformation theory for a category of modules over the Chas-Sullivan Lie bialgebra (see also \cite{CS2}) of an oriented $3$-manifold. But the line of thought will not follow the ideas of classical quantization as in \cite{CFP}. 
 
Here is the main idea of the paper: For $M$ a $3$-dimensional manifold, the Chas-Sullivan structure \textit{measures} the oriented intersections between the loops in a $1$-dimensional family with those in a $0$-dimensional familiy (which is just a formal linear combination of free homotopy classes of loops). The resulting pairing takes values in formal linear combinations of free homotopy classes. It is easy to see that the construction extends to \textit{collections}
of loops. It is our goal to study possible refinements of that resulting extended structure. More precisely we  ask the following question:
 \textit{What happens if we try to replace 
collections of free homotopy classes of loops in $M$ by isotopy classes of oriented links in $M$?} 
Note that in physics terms the collections of loops up to homotopy 
represent classical observables while the isotopy classes of
links are the quantum observables.  

Our main method is transversality of families, which is at the heart of Vassiliev theory \cite{V}.
We will show that the process of replacing homotopy classes of maps by equivalence classes of
\textit{transverse objects} (oriented links in the case of $0$-dimensional homology of the mapping space)
naturally deforms the Chas-Sullivan type \textit{intersection theory} into well-known structures in quantum topology.
Thus our approach provides a weakened version of quantization in the category of oriented $3$-manifolds in comparison to the deep results for
cylinders over oriented surfaces due to Goldman \cite{G} and Turaev \cite{T}. Because of the lack of a geometric product structure in an arbitrary oriented $3$-manifold such an extension will not be based on the deformation of algebra structures but the structures of modules over the Chas-Sullivan Lie
algebra.  

In order to define a \textit{natural} refinement of the Chas-Sullivan structure, isotopy has to be weakened to link homotopy of oriented links in $M$. Recall that in a deformation through link homotopy arbitrary self-crossings are allowed but the different components have to be disjoint during a deformation. 
This equivalence relation has first been considered by Milnor in 1952 \cite{M} (also see the recent results by Koschorke \cite{Ko}). 
But it will turn out that we need additional relations
naturally resulting from the geometry of the problem. These relations can be considered as universal versions of integrability conditions 
previously considered in Vassiliev theory \cite{KL}, \cite{L}. Formally this defines a \textit{universal} quotient of the free abelian group generated by the set of link homotopy classes of oriented links in $M$ by the integrability relations. The relations are parametrized by certain immersions with two
double points of dinstinct components. The resulting group turns out be closely related to a well-known skein module, first studied by Hoste and Przytyki in 1990 \cite{HP1}, later by the author in \cite{K1}. We will show that the difference between the universal quotient and the Hoste-Przytycki module is subtle and determined by the chord diagrams of integrability relations.  

In section 1 we state our main results.
In section 2 we review the Chas-Sullivan string topology 
and set up the basic definitions of our theory. 
In section 3 some results of \cite{K1} are reformulated
and new results about the structure of Hoste-Przytycki skein modules are proven. In section 4 the proof of our main theorem 3 (see section 1) is given. In section 5 we study the relation between the universal
string topology deformation and 
the Hoste-Przytycki skein module. 

\noindent \textbf{Acknowledgement}\qua It is a pleasure to thank Charles
Frohman for bringing the work of Moira Chas and Dennis Sullivan to my 
intention. Also, I would like to thank Jozef Przytycki and Dennis
Sullivan for giving me possibilities to talk about some of the ideas of this paper
in their seminars and workshops.

\section{Statement of the main results} 

Throughout let $R$ be a commutative ring with $1$ and let
$\varepsilon :R\rightarrow \mathbb{Z}$ be an epimorphism of commutative rings with $1$. 
An $R$-Lie algebra is a (not necessarily finitely generated) $R$-module with a
Lie-bracket (which is an $R$-bilinear operation satisfying 
Jacobi-identity and antisymmetry). 
 
Let $M$ be an oriented $3$-manifold and let
$map(S^1,M)$ be the space of continuous maps from the circle
$S^1$ into $M$ (the free loop space of
$M$) with the compact open topology.
There is the circle symmetry of rotating the domain. The 
corresponding equivariant homology
$$\mathcal{H}_*(M;R):=H_*^{S^1}(map(S^1,M);R)$$
is isomorphic to the homology of the quotient space of $map(S^1,M)$ by the
$S^1$-action (string homology). In this paper only the pairings (string interactions)
$$\gamma_1: \mathcal{H}_1(M;R)\otimes \mathcal{H}_0(M;R)\rightarrow \mathcal{H}_0(M;R)$$
and 
$$\quad \gamma _L: \mathcal{H}_1(M;R)\otimes \mathcal{H}_1(M;R)\rightarrow
\mathcal{H}_1(M;R)$$
are needed.
Here $\gamma _L$ is a Lie-bracket on the $1$-dimensional
equivariant homology and $\gamma_1$ equips the $R$-module $\mathcal{H}_0(M;R)$ with the structure of a module over the $R$-Lie algebra $\mathcal{H}_1(M;R)$.
We let $\gamma :=\gamma_1$ denote the pairing in the case $R=\mathbb{Z}$.  

The module $\mathcal{H}_0(M;R)$ is free with basis the set of conjugacy classes of $\pi_1(M)$. It follows from homotopy theory and the results of \cite{CS1} that 
for an oriented $3$-manifold with infinite fundamental group and contractible universal covering, $\gamma_1$ and $\gamma_L$ are the only possibly non-trivial string topology pairings. For $M$ a hyperbolic $3$-manifold both pairings are trivial while the string topology pairing $\gamma $ is known to be non-trivial for most Seifert-fibred $3$-manifolds (compare \cite{CS1}, and \cite{K1}) for the case of a $3$-manifold with boundary).  

The pairing $\gamma $ describes  \textit{interactions} of a given 
loop with the loops in a $1$-dimensional family. This easily extends to 
a pairing defined for a \textit{collection of loops} (in particular oriented links) 
and a $1$-dimensional family. 

Let $\mathcal{L}(M)$ be the set of link homotopy classes of oriented links in $M$ (including the empty link $\emptyset$). 
In section 2 we will define for each $3$-manifold $M$
the notion of a \textit{geometric $R$-deformation} of $\gamma $. 
This is roughly a quotient module $A$ of $R\mathcal{L}(M)$ 
equipped with a pairing $\gamma_A :\mathcal{H}_1(M;R)\otimes A\rightarrow A$ deforming $\gamma $ in a specific 
way based on \textit{transversality}. For a formal definition compare 
definition 3 in section 2.   
It will be an immediate consequence of definition 2 in section 2 that 
$\gamma_A$ induces on $A$ the structure of an $R$-Lie algebra over 
$\mathcal{H}_1(M;R)$.     

Homotopy skein modules are natural quotients of $R\mathcal{L}(M)$ (see \cite{HP1}, \cite{P1}, \cite{K1}) for 
suitable rings $R$. Let $R=\mathbb{Z}[z]$. Then the Hoste-Przytycki module   
$\mathcal{C}(M)$ is defined as the quotient of $R\mathcal{L}(M)$ by the
submodule generated by all linear combinations 
$K_+-K_--zK_0$. Here $K_0$ is the oriented smoothing of two links
$K_{\pm}$ and the three links differ only in the interior of some oriented $3$-ball in
$M$.  

\begin{center}
\epsfig{file=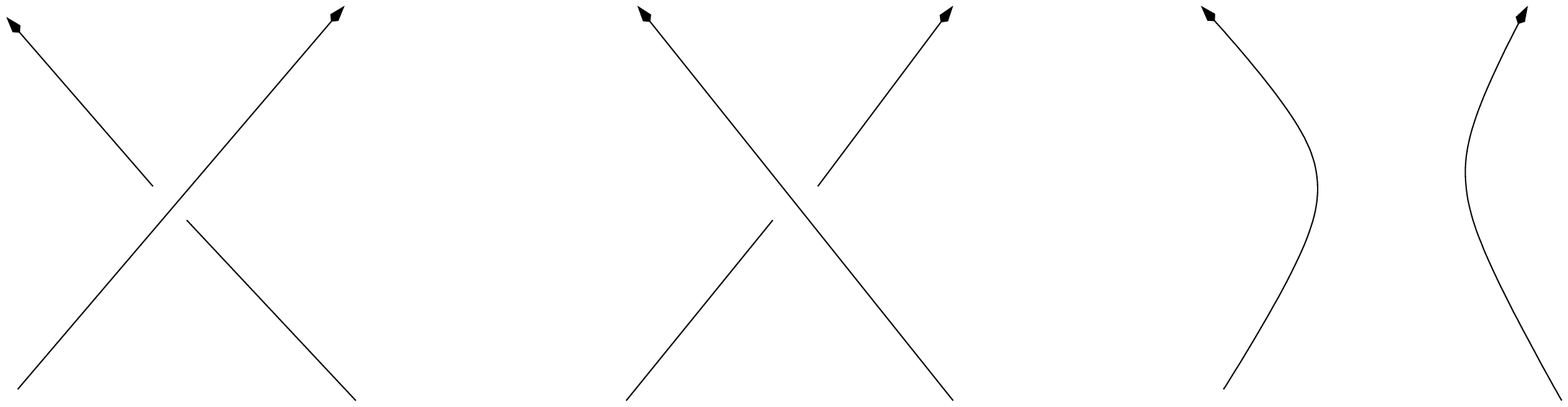,height=0.6in,width=3.2in}
\end{center}
{\small\qquad \qquad \qquad \qquad\quad $K_+$ \ \quad \qquad \qquad \qquad $K_-$ \qua \quad \qquad \qquad \quad $K_0$}

Moreover, the \textit{crossing arcs of} $K_{\pm}$ in the ball belong to \textit{distinct} 
components. The skein relations above are called \textit{homotopy Conway relations}.

\begin{theorem}
The Hoste-Przytycki skein module is a geometric $R$-deformation of $\gamma $ for $R=\mathbb{Z}[z]\rightarrow \mathbb{Z}$ defined by $z\mapsto 0$. In particular there is a pairing 
$$\gamma_{\infty }: \mathcal{H}_1(M;R)\otimes \mathcal{C}(M)\rightarrow \mathcal{C}(M),$$
which equips $\mathcal{C}(M)$ with the structure of a
module over the $R$-Lie algebra $\mathcal{H}_1(M;R)$.
\end{theorem}

The notation $\gamma_{\infty }$ will become clear in section 3. 
Recall from \cite{K1} (see also \cite{K2}) that $\mathcal{C}(M)$ is free if and
only if $\mathcal{C}(M)$ is torsion free (if and only if $\mathcal{C}(M)$ is
isomorphic to the free module on a natural basis of standard links).   

\begin{theorem} The $R$-module $\mathcal{C}(M)$ is
free if and only if the pairing $\gamma _{\infty }$ is trivial.
\end{theorem}

It will be shown in section 3 that the structure of the skein module is inductively encoded in a sequence of geometric $R$-deformations of $\gamma $. 
In definition 4 of section 2 the notion of a \textit{universal} geometric $R$-deformation of $\gamma $ will be defined in the obvious way. For a given $3$-manifold $M$ the universal geometric $R$-deformation is unique up to isomorphism of $R$-Lie algebras over $\mathcal{H}_1(M;R)$.    
The construction is very similar to skein modules. We give the description at first for 
$\varepsilon =id: R=\mathbb{Z}\rightarrow \mathbb{Z}$.  

Consider the set of immersions of circles in $M$ with intersections given by two double points of distinct components.
We will assume  that the unit vectors tangent to the two branches at a double point span a plane in the corresponding tangent space of $M$. By abuse of notation we call these immersions \textit{transversal}. Two transversal immersions are equivalent if they can be deformed into each other by ambient isotopy of the $3$-manifold and smooth homotopies allowing arbitrary self-intersections of components in the complement of the two double points. For $i=2,3,4$ 
let $\mathcal{K}^{i}(M)$ denote the set of equivalence classes of transversal
immersions as above with the double points occuring in $i$ components.
(Equivalently this can be expressed by saying that the two chords in the chord diagrams determined by some immersion in $\mathcal{K}^{i}(M)$ have end points 
on $i$ components).
For $K_{**}\in \mathcal{K}^{i}(M)$ let $(**)$ indicate the two double points and let
$K_{0\pm},K_{\pm 0}$ denote the four natural resolutions of
the two double-points. Let  $K_{+0}-K_{-0}-K_{0+}+K_{0-}\in \mathbb{Z}\mathcal{L}(M)$ be the \textit{integrability element} 
determined by the immersion.
The notation follows \cite{L}, where those relations appear as natural obstructions in the integration of invariants of
finite type in oriented $3$-manifolds. The integrability elements from immersions in $\mathcal{K}^{2}(M)$ are trivial because $K_{0+}=K_{0-}$ and $K_{+0}=K_{-0}$ holds in this case.  

Let $\mathcal{W}(M;\mathbb{Z})$ be the quotient of $\mathbb{Z}\mathcal{L}(M)$ by the subgroup generated by integrability elements from all transverse immersions in 
$\mathcal{K}^{3}(M)$.  Then $\mathcal{W}(M;\mathbb{Z})$ is called the \textit{Chas-Sullivan module} because of the following theorem. Here, for the formal
definition of \textit{universal} geometric $\mathbb{Z}$-deformation, compare definition 4 in section 2.  

\begin{theorem}
The abelian group $\mathcal{W}(M;\mathbb{Z})$ together with a natural pairing
$$\gamma _u : \mathcal{H}_1(M;\mathbb{Z})\otimes \mathcal{W}(M;\mathbb{Z})\rightarrow \mathcal{W}(M;\mathbb{Z}).$$
is the universal geometric $\mathbb{Z}$-deformation of $\gamma $. 
Moreover, $\mathcal{W}(M;R):=\mathcal{W}(M)\otimes R$ together with
the obvious tensor product pairing $\gamma_u^R$ (defined by $R$-linear extension) is the universal
geometric $R$-deformation of $\gamma $.
\end{theorem}

\begin{remark}
(a)\qua Let $D\cong R\mathcal{L}(M)/d \cong \mathcal{W}(M;R)/d'$
be a quotient of the universal geometric $R$-deformation (equivalently the natural
projection $R \mathcal{L}(M) \rightarrow D$ factors through $\mathcal{W}(M;R)$).
Then $\gamma_u^R$ induces on $D$ the structure of a geometric $R$-deformation of $\gamma $ 
if and only if 
the \textit{locality condition}
$$\gamma _u^R(\mathcal{H}_1(M;R)\otimes d')\subset d'.$$
holds.
This will be satisfied e.g.\ for all skein module quotients of
$\mathcal{W}(M;R)$ because the defining skein relations can be assumed to take place in balls separated from intersections with singular tori.

\noindent (b)\qua Assume that $M$ is a $3$-manifold such that every essential (i.e.\ $\pi_1$-injective) torus map is homotopic into $\partial M$. Then it follows from the arguments in \cite{K1} and the construction of $\mathcal{W}(M;\mathbb{Z})$ that $\gamma _u$ is trivial (compare section 4).
  
\noindent (c)\qua It is possible to extend $\gamma_L$ to a natural $R$-Lie algebra structure on the first homology module of the quotient of $\cup_{r\geq 0}map(\cup_rS^1,M)$ by $S^1$-actions 
and permutations of components. The universal module with respect to the structure
of a module over \textit{this} $R$-Lie algebra is the module $\mathcal{W}^+(M;R)$ defined as the quotient of $R\mathcal{L}(M)$ by all possible integrability relations (i.e.\ from immersions in $\mathcal{K}^3(M)\cup \mathcal{K}^4(M)$). 
\end{remark}
  
\begin{example} The Hoste-Przytycki module $\mathcal{C}(M)$
is a quotient of $\mathcal{W}(M;R)$ satisfying the locality condition above. In fact, 
$K_{+0}-K_{-0}=zK_{00}=K_{0+}-K_{0-}\in \mathcal{C}(M)$
holds (\textit{also} for all transverse immersions). Thus all integrability elements map trivially
 into $\mathcal{C}(M)$. It follows that the natural projection 
$R \mathcal{L}(M)\rightarrow \mathcal{C}(M)$ factors through $\mathcal{W}(M;R)$.
In particular theorem 3 implies theorem 1. 
\end{example}

\begin{example}
Not all homotopy skein modules are geometric $R$-deformations of $\gamma $: Let $\varepsilon : R=\mathbb{Z}[q^{\pm},z]\rightarrow \mathbb{Z}$ be defined by
$z\mapsto 0, q\mapsto 1$. The $q$-homotopy skein module (compare \cite{P1}) is defined as the quotient of $R\mathcal{L}(M)$ by the relations 
$q^{-1}K_+-qK_-=zK_0$ for crossings of distinct components. 
In this case not all integrability elements  are contained in the kernel of the
projection of the free module onto the skein module.
Thus the $q$-deformed module is \textit{not} a quotient of $\mathcal{W}(M;R)$. In fact the element
$K_{+0}-K_{-0}-K_{0+}+K_{0-}$ maps to $(q^2-1)(K_{-0}-K_{0-})$ (for $K_{**}\in \mathcal{K}^{3}(M)$). This is a non-trivial element of the skein module in general (see \cite{K1} and \cite{P1}). Note that according to \cite{K1} there are relations in the skein module resulting from $\mathcal{K}^{2}(M)$, which imply 
the existence of torsion in the $q$-homotopy skein module if $\pi_1(M)$ is not abelian (for details see \cite{K1}). 
\end{example}

For each commutative ring $R$ with $1$, the identy element defines the 
ring inclusion $\mu :\mathbb{Z}\rightarrow R$. It induces homomorphisms of abelian groups $\mathbb{Z}\mathfrak{S}\rightarrow R\mathfrak{S}$ for each set $\mathfrak{S}$. Also there are homomorphisms $B\rightarrow B\otimes_{\mathbb{Z}}R$ for each abelian group $B$ defined by $x\mapsto x\otimes 1$. This induces homomorphisms  
of abelian groups 
$$\mathcal{H}_1(M;\mathbb{Z})\rightarrow \mathcal{H}_1(M;\mathbb{Z})\otimes_{\mathbb{Z}}R\cong
\mathcal{H}_1(M;R)$$ 
using  
the universal coefficient theorem, and 
$$\mathcal{W}(M;\mathbb{Z}) \rightarrow \mathcal{W}(M;\mathbb{Z})\otimes_{\mathbb{Z}}R\cong \mathcal{W}(M;R).$$ 
All the induced homomorphisms are denoted $\mu_*$.
It follows from theorem 3 and the definition of universal geometric $R$-deformation as given in section 2, definition 3, that there is the commuting diagram:
$$
\begin{CD}
\mathcal{H}_1(M;\mathbb{Z})\otimes \mathcal{W}(M;\mathbb{Z})@>\gamma _u >> \mathcal{W}(M;\mathbb{Z})\\
@V{\mu_*\otimes \mu_*}VV @V{\mu_*}VV \\
\mathcal{H}_1(M;R)\otimes \mathcal{W}(M;R)@>\gamma_u^R>> \mathcal{W}(M;R)
\end{CD}
$$  

Let $\mathcal{L}_{\centerdot}(M):=\mathcal{L}(M)\setminus \{\emptyset \}$.
All the constructions of modules, homomorphisms and pairings (given in detail in section 2) can be performed by replacing $\mathcal{L}(M)$ by $\mathcal{L}_{\centerdot}(M)$. For each such module $A$ let $A_{\centerdot}$ be the resulting \textit{reduced} module. Similar notation applies to homomorphisms and pairings.
E.g.\ recall from \cite{HP1} that $\mathcal{C}(M)=\mathcal{C}_{\centerdot}(M)\oplus R$ where $R$ is the submodule generated by the empty link and $\mathcal{C}_{\centerdot}(M)$ is the \textit{reduced Hoste-Przytycki module}. 

Because of universality, each geometric $R$-deformation $A$ of $\gamma $ comes equipped with a unique 
$R$-epimorphism $\mathcal{W}(M;R)\rightarrow A$ and a corresponding epimorphism of reduced modules $\mathcal{W}_{\centerdot}(M;R)\rightarrow A_{\centerdot}$. 

\begin{definition}
A geometric $R$-deformation $A$ (and also the reduced $R$-defor-\
mation $A_{\centerdot}$) of $\gamma $
with pairing $\gamma_A$ is called \textit{integral} if the \textit{integrality homomorphism}
$$\mu _{A_{\centerdot}}: \mathcal{W}_{\centerdot}(M;\mathbb{Z})\rightarrow \mathcal{W}_{\centerdot}(M;R)\rightarrow A_{\centerdot}$$
is surjective.  
\end{definition}

\begin{remark} 
(a)\qua Because the homomorphism $\mathbb{Z}\mathcal{L}_{\centerdot}(M)\rightarrow \mathcal{W}_{\centerdot}(M;\mathbb{Z})$ is onto it follows easily that 
a geometric $R$-deformation $A$ is integral if and only if the composition
$$\mathbb{Z}\mathcal{L}_{\centerdot}(M)\rightarrow R\mathcal{L}_{\centerdot}(M)\rightarrow {A}_{\centerdot}$$
is onto.    

\noindent (b) For an integral geometric $R$-deformation $A$, it follows from the commutative diagram above that the pairing
$\gamma_{A_{\centerdot}}$ is defined by $R$-linear extension in the first factor
from the geometric $\mathbb{Z}$-deformation given by $\mu _{A_{\centerdot}}$ 
(and remark 1 (a) above) for the abelian group $A_{\centerdot}$.    
For a non-trivial ring extension $R\rightarrow \mathbb{Z}$ this endows a reduced integral geometric $R$-deformation with an $R$-module structure, which is  \textit{compatible} with the induced structure of a geometric $\mathbb{Z}$-deformation of $\gamma$. 
\end{remark}

By applying the skein relation $zK_0=K_+-K_-$ to a connected sum of a link $K$
and a Hopf link (with $K_0=K$) it follows from remark 2 (a) that $\mathcal{C}(M)$ is integral. In section 5 the following result is proved:

\begin{theorem}
Suppose that $M$ is a submanifold of a rational homology $3$-sphere. If $H_1(M)\neq 0$  
then the integrality epimorphism
$$\mathcal{W}_{\centerdot}(M;\mathbb{Z})\rightarrow \mathcal{C}_{\centerdot}(M)$$
is not injective.
\end{theorem} 

Recall that $\mathcal{W}^+(M;\mathbb{Z})$ be the quotient of $\mathbb{Z}\mathcal{L}(M)$ by \textit{all} possible integrability relations.
In section 5 we will prove the isomorphism
$\mathcal{W}_{\centerdot}^+(M;\mathbb{Z})\cong \mathcal{C}_{\centerdot}(M)$. 
Combined with remark 1 (c) from above and remark 3 from section 4 it follows that the reduced Hoste-Przytycki module is the reduced
universal geometric $\mathbb{Z}$-deformation of $\gamma $ (as module over the extended Lie algebra as described in remark 1 (c) ). 

As a by-product of the proof of theorem 2 we construct the first example of $\mathbb{Z}$-torsion in $\mathcal{C}(M)$.
For general discussions of torsion in skein modules see \cite{HP2}, \cite{P2} and \cite{P3}.  

\begin{ztor} 
Let $K$ be the $2$-component link in 
$S^2\times S^1$ defined from two oriented parallel copies of $\{*\}\times S^1$. Then 
in $\mathcal{C}(S^2\times S^1)$ the relations 
$zK\neq 0$ but $2(zK)=0$ hold.   
\end{ztor}

\begin{conjecture}
(a)\qua For each oriented $3$-manifold, the epimorphism $\mu_{\mathcal{C}_{\centerdot}(M)}: \mathcal{W}_{\centerdot}(M;\mathbb{Z})\rightarrow \mathcal{C}_{\centerdot}(M)$ is not an isomorphism.

\noindent (b)\qua Assume that each essential torus map into $M$ is homotopic into $\partial M$. Then $\mathcal{W}(M;\mathbb{Z})$ is free and the universal
pairing is trivial.     

\noindent (c)\qua For a given $3$-manifold $M$, if $\gamma $ is trivial then each geometric $R$-deformation of $\gamma $ is trivial. (This is of course an important problem in the understanding of the structure of $\mathcal{C}(M)$.) 
\end{conjecture}

All constructions in this paper are functorial in the obvious way with 
respect to embeddings of oriented $3$-manifolds. In
theorem 3 also obvious functoriality with respect to homomorphisms of
commutative rings with $1$ holds.

\section{Geometric deformations of the Chas-Sullivan pairing}

We first review the necessary background from string topology. Compare \cite{CS1}, in particular section 1 and section 6. 
 
For $R$ a commutative ring with $1$ and $A$ some $R$-module let
$SA$ denote the symmetric algebra.
 
Let $\mathbb{H}_*(M;R):=H_*(map(S^1,M);R)$ be
the usual homology of the free loop space of $M$ (loop homology).
By the universal coefficient theorem 
$H_i(\ - \ ;R)\cong H_i(\ - \ )\otimes_{\mathbb{Z}} R$
for $i=0,1$, and this holds both for string and loop homology. 

We describe the string homology pairings over $\mathbb{Z}$, the pairings over $R$ 
then are defined by linear extension. 

The quotient mapping
$map(S^1,M)\rightarrow map(S^1,M)/S^1$, which restricts to an $S^1$-fibration over the subspace of non-constant loops, induces the exact Gysin-sequence:
$$\rightarrow \mathbb{H}_1(M;\mathbb{Z})\rightarrow \mathcal{H}_1(M;\mathbb{Z})\rightarrow
0\rightarrow \mathbb{H}_0(M;\mathbb{Z})\rightarrow \mathcal{H}_0(M;\mathbb{Z})\rightarrow 0.$$
Thus the elements of $\mathcal{H}_0(M;\mathbb{Z})$ can be identified with elements of 
$\mathbb{H}_0(M;\mathbb{Z})$, and the elements of $\mathcal{H}_1(M;\mathbb{Z})$  can be lifted to $\mathbb{H}_1(M;\mathbb{Z})$.
So each element of $\mathbb{H}_0(M;\mathbb{Z})$ is
represented by a linear combination of oriented knots in $M$ (
$0$-dimensional family). Correspondingly elements of $\mathbb{H}_1(M;\mathbb{Z})$ are represented by linear combinations of $S^1$-families of mappings from a circle into $M$.
The image of each such an $S^1$-family in $M$ is a singular torus. The corresponding maps can be assumed in transversal position with respect to knots representing elements of 
$\mathcal{H}_0(M;\mathbb{Z})$ and singular tori respresenting other elements of $\mathcal{H}_1(M;\mathbb{Z})$. Then at transversal intersection points new loops are defined by going around the first loop and then around the second one (usual loop composition or loop multiplication), with the loops determined by the two preimages of the intersection point.
(The parametrization of the singular torus as a family of loops in the $3$-manifold is used here in some essential way.)
This defines new $0$-dimensional and $1$-dimensional families of maps. The construction can be carried out for each pair of $S^1$-family and oriented knot, and for each pair of $S^1$-families. The resulting linear combinations are used to define $\gamma $ and $\gamma_L$. 

Next we replace the second factor in the domain of $\gamma $ by the $0$-dimensional homology of the space of maps of circles into $M$ with arbitrary numbers of components.   
Note that $\mathcal{H}_0(M;\mathbb{Z})\cong \mathbb{H}_0(M;\mathbb{Z})\cong \mathbb{Z}\hat{\pi}(M)$, where
$\hat{\pi}(M)$ is the set of conjugacy classes of elements in the
fundamental group $\pi_1(M)$. By the isomorphism above $\gamma $ is a pairing 
$ \mathcal{H}_1(M;\mathbb{Z})\otimes \mathbb{Z}\hat{\pi }(M)\rightarrow \mathbb{Z}\hat{\pi}(M)$.

Let $Map(S^1,M):= \cup_{r\geq 0}map(\cup_rS^1,M)/\Sigma_r$, where
the permutation group $\Sigma_r$ acts by permuting the components of the domain. Then 
$H_0(Map(S^1,M);\mathbb{Z})\cong S\mathbb{Z}\hat{\pi}(M)$
and the pairing $\gamma $ naturally extends to the pairing (still denoted)
$$\gamma  : \mathcal{H}_1(M;\mathbb{Z})\otimes S\mathbb{Z}\hat{\pi}(M)\rightarrow S\mathbb{Z}\hat{\pi}(M)$$
using the \textit{Poisson identity}
$$\gamma (a\otimes xy)=\gamma (a\otimes x)y+x\gamma (a\otimes y).$$
It follows  by some easy computation that $\gamma $ induces on $S\mathbb{Z}\hat{\pi}(M)$ the structure of a module over the Lie algebra $\mathcal{H}_1(M;\mathbb{Z})$ (compare also lemma 2 below). 
By definition the extension is the linear combination of the loop compositions at all intersections of a parametrized singular torus with the
components of a \textit{link} in $M$ (whose homotopy classes correspond to a given element in $S\mathbb{Z}\hat{\pi}(M)$). 
The two curves, which are composed at some intersection point of a singular torus with the link, are the corresponding curve in the family $S^1\rightarrow map(S^1,M)$ and the corresponding component of the link. Note that a given loop of the $S^1$-family can have transverse intersections with several components of the link. Then the collection of loops resulting from loop composition for one of those intersection points is a map of circles whose images in $M$ are \textit{not} disjoint.  
It is the starting observation of this paper that this can be avoided using \textit{higher order transversality} following \cite{L}.   

\begin{definition} 
Let $U\subset Map(S^1,M)$ be an open subspace.
Assume that a $0$-chain in $Map(S^1,M)$ is represented by elements in $U$
and a $1$-chain in $map(S^1,M)\subset Map(S^1,M)$ is given. Then assume that the $0$-chain and the $1$-chain
(more precisely the representing families in $Map(S^1,M)$) are approximated such that
their intersection is transversal in the sense of Chas-Sullivan (i.e.\ the 
corresponding obvious evaluations are transverse maps to $M$) and the following two 
additional conditions hold:  
(i) The $0$-chain is approximated by maps in $U$ and (ii) all maps of  
the $0$-dimensional family, which result by smoothing or loop composition at intersection points according to the definition given in \cite{CS1}, are contained in $U$. Then we say that the resulting $0$-dimensional family (and its homology class) is defined by the \textit{Chas-Sullivan procedure} (with respect to $U$). 
\end{definition}

\begin{remark}
Loop composition is the standard operation in homotopy theory while 
smoothing is the standard operation in link theory.  
There is no difference in the definition of $\gamma $ because the resulting $0$-dimensional families defined by the Chas-Sullivan procedure are the same (the corresponding links are link homotopic). It can easily be seen that $\gamma_L$ can also be defined using smoothing.  
But the \textit{general} string \textit{interactions} in \cite{CS1} a priori require the use of loop composition since the resolutions take place along high dimensional cells of intersections.
\end{remark}

The main example considered in this paper is $U=\mathcal{D}(S^1,M)$, the space of maps with disjoint images of its components, which is denoted the \textit{Milnor space} of $M$. The path-components of this space are the link homotopy classes of oriented links in $M$ \cite{M}. 
Then $H_0(\mathcal{D}(S^1,M))\cong \mathbb{Z}\mathcal L(M)$. 
The inclusion 
$\mathcal{D}(S^1,M)\subset Map(S^1,M)$ induces the homomorphism $$\mathfrak{h}: \mathbb{Z}\mathcal{L}(M)\rightarrow S\mathbb{Z}\hat{\pi}(M),$$
which assigns to each link the unordered sequence of the homotopy classes of its components.

In order to achieve the transversality necessary for the Chas-Sullivan procedure
we define an operation, which takes 
the union of two families of maps of circles in $M$ with the same parameter space. Let 
$\widetilde{Map}(S^1,M):=\cup_{r\geq 0} map(\cup_rS^1,M)$ with the natural projection (a covering map away from maps with two equal component maps)
$\widetilde{Map}(S^1,M)\rightarrow Map(S^1,M)$.
Next consider the map
$$\widetilde{Map}(S^1,M)\times \widetilde{Map}(S^1,M)\rightarrow \widetilde{Map}(S^1,M),$$
defined by the disjoint union of maps
$$map(\cup_rS^1,M)\times map(\cup_sS^1,M)\rightarrow map(\cup_{r+s}S^1,M).$$

\begin{remark}
There are the induced maps after taking the quotient by permutation actions:
$$map(\cup_rS^1,M)/\Sigma_r\times map(\cup_sS^1,M)/\Sigma_s\rightarrow map(\cup_{r+s}S^1,M)/\Sigma_{r+s},$$
and thus $Map(S^1,M)\times Map(S^1,M)\rightarrow Map(S^1,M)$. 
The induced map in homology and the K\"{u}nneth homomorphism define the 
homomorphism of degree $0$ (graded commutative multiplication):
$$H_*(Map(S^1,M);\mathbb{Z})\otimes H_*(Map(S^1,M);\mathbb{Z})\rightarrow H_*(Map(S^1,M);\mathbb{Z}),$$
which restricts to the standard multiplication in $0$-dimensional homology
$$S\mathbb{Z}\hat{\pi}(M)\otimes S\mathbb{Z}\hat{\pi}(M)\rightarrow  
S\mathbb{Z}\hat{\pi}(M).$$
used above. (Note that $S\mathbb{Z}\hat{\pi}(M)$ is multiplicatively generated by $\mathbb{Z}\hat{\pi}(M)$.) This multiplication in homology and similar structures  
will be essential for developing generalizations of
the deformation theory discussed in future work.  
\end{remark}

Now given an element in $\widetilde{Map}(S^1,M)$ represented by some embedding $e$, consider the constant loop $S^1\rightarrow \widetilde{Map}(S^1,M)$ in $e$.
Its union with a $1$-dimensional family $\ell : S^1\rightarrow map(S^1,M)\subset \widetilde{Map}(S^1,M)$ defines a new $1$-dimensional family $\ell_e$. The following result is easily proved.

\begin{lemma} After approximation of $e$ and $\ell $ it can be assumed that $\ell_e$ is transversal, i.e.\ the maps $\cup_rS^1\rightarrow M$ in the family given by
$\ell_e$ are all embeddings except for a finite number of parameters. For those parameters the maps are immersions with a single transverse double-point.
\end{lemma}

Here we will only need the \textit{weak transversality} version of the previous
transversality lemma for link homotopy (compare \cite{K1}).
Then $\ell_e: S^1\rightarrow Map(S^1,M)$ is
\textit{transversal} as follows: All maps in the family are embeddings
(in particular have disjoint images of components) except 
for a finite set of \textit{singular} parameters. 
At a singular parameter the corresponding map $\cup_rS^1\rightarrow M$ is an immersion (actually the image of an immersion in $\widetilde{Map}(S^1,M)$ projected to $Map(S^1,M)$) with a single \textit{transverse} double-point of distinct components.
Then the linear combination, defined by smoothing at all intersection points, is defined by the \textit{Chas-Sullivan procedure} according to definition 1 and $U=\mathcal{D}(S^1,M)$.  

The following algebraic language is very useful.
Let $A,B$ be modules over commutative rings $R,R'$ and let $\delta :R\rightarrow R'$ be a ring homomorphism. Then a homomorphism $h: A \rightarrow B$ is a homomorphism of abelian groups satisfying
$h(rx)=\delta (r)h(x)$ for all $r\in R$ and $x\in A$. 
If $R$ a commutative ring with $1$ and $\delta : \mathbb{Z}\rightarrow R$ is the injection defined by the unit then a homomorphism is just a homomorphism of underlying abelian groups. 

Now recall that the ring $R$ is equipped with the ring epimorphism
$\varepsilon :R\rightarrow \mathbb{Z}$. There is the induced homomorphism
$$
\begin{CD}
\varepsilon _*: \mathcal{H}_1(M;R) \cong \mathcal{H}_1(M;\mathbb{Z})\otimes_{\mathbb{Z}}R
@>{id \otimes \varepsilon }>>\mathcal{H}_1(M;\mathbb{Z}) \otimes_{\mathbb{Z}} \mathbb{Z} \cong \mathcal{H}_1(M;\mathbb{Z}).     
\end{CD}
$$
In general, homomorphisms from free $R$-modules to free abelian groups induced by $\varepsilon $ will be denoted $\varepsilon_*$. 

\begin{definition}
Consider a quadruple $(A,\alpha_A,\beta_A,\gamma_A)$ with $A$ an $R$-module with epimorphisms $\alpha_A,\beta_A$ as follows:
$$
\begin{CD} 
R\mathcal{L}(M)@>\alpha_A>>A@>\beta_A >>S\mathbb{Z}\hat{\pi}(M),
\end{CD}
$$
such that $\beta_A \circ \alpha_A$ is the epimorphism
$$
\begin{CD}
R\mathcal{L}(M)@>\varepsilon_* >> \mathbb{Z}\mathcal{L}(M)@>{\mathfrak{h}}>> S\mathbb{Z}\hat{\pi}(M),
\end{CD}
$$
and $\gamma_A : \mathcal{H}_1(M;R)\otimes A\rightarrow A$ is $R$-linear and defined 
on representative links by the Chas Sullivan procedure.
Then $(A,\alpha_A,\beta_A,\gamma_A)$, or briefly $\gamma_A$, is called a \textit{geometric $R$-deformation of $\gamma $}.
\end{definition}

\begin{remark}
It follows that there is the commuting diagram
$$
\begin{CD}
\mathcal{H}_1(M;R)\otimes A@>\gamma_A >> A\\
@V{\varepsilon_* \otimes \beta_A}VV @V{\beta_A}VV \\
\mathcal{H}_1(M;\mathbb{Z})\otimes S\mathbb{Z}\hat{\pi}(M)@>\gamma >> S\mathbb{Z}\hat{\pi}(M)
\end{CD}
$$  
\end{remark}

\begin{example} 
The trivial example of a geometric $R$-deformation of $\gamma $ is given by $$A:=H_0(Map(S^1,M);R)\cong R\mathcal{L}(M)/a,$$ where $a$ is the submodule generated by $K_+-K_-$ for all $K\in \mathcal{L}(M)$. Note that there is the canonical homomorphism
$$\beta_A: SR\hat{\pi}(M)\cong H_0(Map(S^1,M);R)\rightarrow H_0(Map(S^1,M),\mathbb{Z})\cong S\mathbb{Z}\hat{\pi}(M)$$
induced by $\varepsilon $. The pairing $\gamma_A$ is defined by $R$-linear extension of $\gamma $. Note that this trivial geometric $R$-deformation is not co-universal in the sense that each geometric $R$-deformation maps onto it.
In fact, the Hoste-Przytycki homotopy skein module maps onto the trivial geometric $R$-deformation if and only if it is free.    
\end{example}

Recall that an $R$-module $A$ is a \textit{module over the $R$-Lie algebra $\mathcal{H}_1(M;R)$}
if there is defined an $R$-homomorphism:
$\mathcal{H}_1(M;R)\rightarrow End_R(A)$,
which maps the given Lie-bracket $\gamma_L$ on $\mathcal{H}_1(M;R)$ to the endomorphism Lie-bracket. 

The following result is immediate from the argument on page 2 of \cite{CS1} proving the bracket property in string topology, and the fact that the deformation is geometric.

\begin{lemma}
If $A$ is a geometric $R$-deformation of $\gamma $ then $\gamma_A$ 
induces on $A$ the structure of an $R$-Lie algebra over $\mathcal{H}_1(M;R)$.
\end{lemma}

The induced structure on $A$ of a module over the Chas-Sullivan Lie algebra is called a geometric 
$R$-deformation of the module $\mathcal{H}_0(M;\mathbb{Z})$ over the Lie algebra $\mathcal{H}_1(M;\mathbb{Z})$.   

The result of lemma 2 is central for the constructions in this paper. Thus 
for the convenience of the reader we give the argument of Chas-Sullivan in the language and orientation conventions used in this paper. The argument also hints into a direction of a general homology and intersection theory of transversal chains in the mapping spaces considered here.
 
\begin{proof}[Proof] For the proof let $\gamma :=\gamma_A$. We want to show that the Jacobi-identity holds in the form 
\begin{equation}
\gamma (\gamma_L(a,b),x)=\gamma (a, \gamma (b,x))-\gamma (b,\gamma (a,x))
\end{equation}
for $x\in A$ and and $a,b\in \mathcal{H}_1(M;R)$. Since $A$ is generated by
link homotopy classes of oriented links in $M$ we can assume that $x$ is 
represented by some oriented $r$-component link $K$ in $M$. Moreover $a,b$ can be represented by two families $S^1\rightarrow map(S^1,M)$.   
Using the union operation above,
$\gamma_L(a,b)$ is represented by smoothing at the \textit{singular points} of the resulting family of maps $f_{ab}: S^1\times S^1\rightarrow map(S^1\cup S^1,M)$.
This is a collection of oriented immersed loops $\ell_{ab}$ in $S^1\times S^1$. Now take the union of $f_{ab}$ with the family $S^1\times S^1\rightarrow map(\cup_rS^1,M)$, which maps constantly, using a parametrization of $K$. By approximation of the original $S^1$-families we can assume that the resulting 
family $S^1\times S^1\rightarrow map(S^1\cup S^1\cup \cup_rS^1,M)$ is transversal in Lin's sense (compare \cite{L}, \cite{K1} and section 4). Because of the fact that the family is constant on the last $r$ components the curves of singular points $\ell_{a,x}$, where the first component intersects a component of $K$ is a union of oriented merdian curves $\{*\}\times S^1$ on $S^1\times S^1$. Similarly the curves $\ell_{b,x}$ form a collection of oriented longitudinal curves $S^1\times \{*\}$ on $S^1\times S^1$. 
Now the singular points corresponding to $f_{ba}$ form the collection of oriented immersed loops $\ell_{ba}$ resulting from $\ell_{ab}$ by application of the switch homeomorphism $S^1\times S^1\rightarrow S^1\times S^1$. 
For any two oriented immersed loops $\ell_1,\ell_2$ on $S^1\times S^1$ intersecting
transversely let $\ell_1\cdot \ell_2\in A$ denote their \textit{oriented intersection in $A$}, which is given by the sum of signed smoothings for all intersection points of $\ell_1$ and $\ell_2$. In fact each intersection point determines an immersion in $\mathcal{K}^{2}(M)\cup \mathcal{K}^{3}(M)$ (compare section 4), which determines a link homotopy class of a smoothed link $\pm K_{00}$ for a suitable $K_{**}$.
If $\ell_{a,x}'$ (respectively $\ell_{b,x}'$) denote the images of $\ell_{a,x}$ (respectively $\ell_{b,x}$) under the switch map then $\ell_{a,x}\cdot \ell_{b,x}=-\ell_{a,x}'\cdot \ell_{b,x}'=\ell_{b,x}'\cdot \ell_{a,x}'\in A$, since switching the order of two curves on a torus changes the intersection number, and the corresponding smoothed links coincide. Here the map into $map(\cup_{r+2}S^1,M)$ on the image of the torus under the switch map is defined by composing the switch map with the map on the domain torus. Similarly $\ell_{ab}\cdot \ell_{a,x}=-\ell_{ba}\cdot \ell_{a,x}'$.
Now consider the following computation of intersections in $A$ of oriented curves on $S^1\times S^1$:
$$
\begin{aligned}
\  & \ell_{ab}\cdot (\ell_{a,x}  \cup \ell_{b,x})\\ 
=  & \ell_{ab}\cdot \ell_{b,x}+\ell_{a,x}\cdot \ell_{b,x} 
+ \ell_{ab}\cdot \ell_{a,x}-\ell_{a,x}\cdot \ell_{b,x} \\
=  & \ell_{ab}\cdot \ell_{b,x}+\ell_{a,x}\cdot \ell_{b,x}
- \ell_{ba}\cdot \ell_{a,x}'-\ell_{b,x}'\cdot \ell_{a,x}'\\
=  & \ell_{ab}\cdot \ell_{b,x}+\ell_{a,x}\cdot \ell_{b,x}
- (\ell_{ba}\cdot \ell_{a,x}'+\ell_{b,x}'\cdot \ell_{a,x}')
\end{aligned}
$$
Now the first line coincides with the left hand side of equation (1) and the last line coincides with the right hand side of (1). 
\end{proof} 

\begin{remark}
(a)\qua Note that the $R$-module structure does not play any role in the proof of the lemma. This is because we actually prove the Jacobi-indentity on the level of the $R$-homomorphism
$$\mathcal{H}_1(M;R)\otimes R\mathcal{L}(M)\rightarrow A,$$
defined by composition with $\alpha_A$ in the second factor, and 
$$\mathcal{H}_1(M;R)\cong \mathcal{H}_1(M;\mathbb{Z})\otimes_{\mathbb{Z}}R$$
holds. (This
requires lifting elements from $A$ to $R\mathcal{L}(M)$.) 
The $R$-module structure itself is only important for 
choosing submodules $a\subset R\mathcal{L}(M)$ for which the image of the pairing above in
$R\mathcal{L}(M)/a\cong A$ is well-defined. 

\noindent (b)\qua The union of all closed oriented curves $\ell_{ab}$, $\ell_{a,x}$ and $\ell_{b,x}$ in the proof of the lemma comes naturally equipped with a transverse map into $Map(S^1,M)$ and thus is a \textit{transverse} $1$-dimensional cycle (see lemma 1 above).
\end{remark}

\begin{definition}
A \textit{universal geometric $R$-deformation}
of $\gamma $  
is a geometric $R$-deformation \newline 
$(C,\alpha_C,\beta_C,\gamma_C)$ such that for each geometric 
$R$-deformation $(A,\alpha_A ,\beta_A ,\gamma_A )$ of $\gamma $ the following holds:
$\alpha_A=\alpha_C\circ \delta_A$ for a unique epimorphism 
$\delta_A: C\rightarrow A$, which fits into the commuting diagram: 
$$
\begin{CD}
\mathcal{H}_1(M;R)\otimes C@>\gamma_C >> C\\
@V{id\otimes \delta_A }VV @V{\delta_A }VV \\
\mathcal{H}_1(M;R)\otimes A@>\gamma_A>> A
\end{CD}
$$
\end{definition}           

It follows easily from the definitions that a universal geometric $R$-deformation of $\gamma $ is \textit{unique up to isomorphism of Lie algebra modules}. 

\begin{remark}
By theorem 3 the universal geometric $R$-deformation is the quotient of $R\mathcal{L}(M)$ by the $R$-submodule generated by integrability relations. Note that the $R$-homomorphism $\mathcal{W}(M;R)\rightarrow SR\hat{\pi}(M)$ defined from universality of $\gamma_u^R$ is the natural one induced by $\mathfrak{h}$ (compare the beginning of section 4).    
\end{remark}

\section{Homotopy skein modules and string topology pairings}

Throughout this section let $R=\mathbb{Z}[z]$.

In \cite{K1} the author has given a presentation of
$\mathcal{C}(M)$ as a quotient of the free module $SR\hat {\pi}(M)$
by a submodule of relations defined geometrically by formal sums 
resulting from smoothing crossings at points of  
intersections of singular tori with \textit{standard links} in $M$
and expanding in terms of standard links. (The standard links are
identified with elements of $SR\hat{\pi}(M)$.) 
For our purposes it is useful to state the presentation 
in a form involving the filtration by the number of components.

Let $r$ be a non-negative integer and let $\mathcal{L}_r(M)$
be the set of link homotopy classes of oriented links in $M$ with
$\leq r$ components. The module $\mathcal{C}_r(M)$ is
defined as the quotient of $R\mathcal{L}_r(M)$ by the submodule $c_r(M)$ of homotopy Conway skein relations for all triples $(K_{\pm},K_0)$ with
$K_{\pm} \in \mathcal{L}_r(M)$. (Note that $K_0\in \mathcal{L}_{r-1}(M)$ for all relations in $c_r(M)$.)

\vskip .1in

There are obvious homomorphisms $\iota_r: \mathcal{C}_{r}(M)\rightarrow
\mathcal{C}_{r+1}(M)$ such that $\mathcal{C}(M)$ is the direct limit
with respect to the $\{\iota_r\}$. In particular there are natural homomorphisms $\mathcal{C}_r(M)\rightarrow \mathcal{C}(M)$ induced by inclusions $\mathcal{L}_r(M)\rightarrow \mathcal{L}(M)$. But the projection $R\mathcal{L}(M)\rightarrow R\mathcal{L}_r(M)$, which is defined by mapping all links with $\geq r+1$ components trivially, does in general not induce a homomorphism from $\mathcal{C}(M)$ to $\mathcal{C}_r(M)$.    
   
For $r\geq 0$ let $\mathfrak{b}_r(M)$ be the set of
unordered sequences of elements of $\hat{\pi }(M)$ of length $r$.
Then $\mathfrak{b}_0(M)$ is the set with the only element given by the empty sequence $\emptyset $ of length $0$. 
Let 
$$\mathfrak{b}(M):=\cup_{r\geq 0}\mathfrak{b}_r(M)$$
be the set of all unordered sequences of elements in $\hat{\pi}(M)$.  
This is the natural basis of monomials of $SR\hat{\pi}(M)$ 
(as $R$-module).
If $\mathfrak{b}_{\centerdot}(M):=\mathfrak{b}(M)\setminus \{\emptyset \}$.
then $R\mathfrak{b}_{\centerdot}(M)$ is also an $R$-algebra but without a unit
element. 

\begin{remark}
If $\mathcal{L}(M)$ is replaced by $\mathcal{L}_{\centerdot}(M)$ and 
$\mathfrak{b}(M)$ is replaced by $\mathfrak{b}_{\centerdot}(M)$ then the definitions and results of section 2 hold for reduced modules, as noticed in the introduction. 
\end{remark}

Choose a
representative oriented link $K_{\alpha }$ for each element $\alpha \in \mathfrak{b}_r(M)$.
The link homotopy classes of these links form a generating set of
$\mathcal{C}_r(M)/\iota_{r-1}(\mathcal{C}_{r-1}(M))$.
The epimorphisms
$$\sigma_r : R(\cup_{0\leq j\leq r}\mathfrak{b}_j(M))\rightarrow \mathcal{C}_r(M)$$
induce the \textit{standard link epimorphism}
$$\sigma : SR\hat{\pi}(M)\rightarrow \mathcal{C}(M).$$

For $\alpha \in \mathfrak{b}_r(M)$ let $f: S^1\rightarrow M$ be a
knot representing some element $\alpha_0 \in \hat{\pi}(M)$, which appears in
$\alpha $.
Let $\delta \in \pi_1(map(S^1,M),f)$ be represented by a map
$h: S^1\rightarrow map(S^1,M)$. Consider a representative link $g$ with
componenents whose homotopy classes correspond to the elements of
$\alpha \setminus \alpha_0$. Then consider 
$S^1\rightarrow map(\cup_{r+1}S^1,M)$, which is defined by taking the
union of $h$ with the constant mapping to $g$, as defined in section 2. 
After approximation of $h$ it can be assumed that the corresponding family
satisfies the transversality conditions discussed in section 2 (see also section 4). Thus for a
finite number of parameters, 
there will be transverse double-points of $h$ with a component of $g$,
and these are the only singularities. 
Then, for each \textit{singular}
parameter, consider the smoothed link, which results by smoothing
the corresponding crossing. Expand this link in terms of standard
links,
which are identified with the corresponding elements of $\mathfrak{b}(M)$.
This defines elements $\Theta (\alpha ,\alpha_0,h)\in SR\hat {\pi}(M)$ which,
after multiplication by $z$, generate the module $\mathcal{R}\subset 
SR\hat{\pi}(M)$ of relations. Note that the relations resulting from
$\alpha \in \mathfrak{b}_r(M)$ 
are well-defined modulo those relations resulting from elements
in $\mathfrak{b}_{r'}(M)$ for $r'<r$.

The ring epimorphism $R \rightarrow \mathbb{Z}$, defined by
$z\mapsto 0$, induces the epimorphism
$$\mathcal{C}(M)\rightarrow S\mathbb{Z}\hat{\pi }(M).$$ 
Its
composition with $\sigma $ is the natural
epimorphism $SR\hat{\pi}(M)\rightarrow S\mathbb{Z}\hat{\pi}(M)$
induced by the the ring epimorphism.

The following is a reformulation of results from \cite{K1}.
It shows that the structure of the Hoste-Przytycki modules is determined 
by pairings, which form a tower of geometric $R$-deformations of $\gamma $.   

\begin{theorem} 
For each non-negative integer $r$
there is a well-defined
pairing
$$\gamma_r: \mathcal{H}_1(M;R)\otimes \mathcal{C}_r(M)\rightarrow
\mathcal{C}_r(M),$$
which is a geometric $R$-deformation of $\gamma $, and there is an isomorphism
$$\mathcal{C}_{r+1}(M)\cong (\mathcal{C}_r(M)/z\cdot  im(\gamma_r))\oplus 
R\mathfrak{b}_{r+1}(M),$$
induced by $\iota_r$ and $\sigma_r$.
\end{theorem}

\begin{proof}[Proof] Let $map_f(S^1,M)$ be the component of $map(S^1,M)$ containing the map $f$ with free homotopy class $[f]\in \hat{\pi}(M)$.
There is the epimorphism 
$$\rho : \bigoplus_{[f]\in \hat{\pi}(M)}\mathbb{Z}\pi_1(map(S^1,M),f)\rightarrow  \bigoplus_{[f]\in \hat{\pi}(M)}H_1(map_f(S^1,M))\rightarrow \mathcal{H}_1(M).$$

The result of theorem 5 follows from the inductive presentation 
of $\mathcal{C}(M)$ in \cite{K1}, using the epimorphism $\rho $ and the fact that
$\Theta (\alpha ,\alpha_0,h)$ only depends on the image of 
$h$ in $\mathcal{H}_1(M)$, which easily follows from the definition (compare also 
the arguments given in section 4, remark 12 (c).    
Alternatively well-definedness of the pairings $\gamma_r$ can be proved from the universal pairing of theorem 3. In fact, the modules $\mathcal{C}_r(M)$ are quotients of $\mathcal{W}(M,R)$: The homomorphism $R\mathcal{L}(M)\rightarrow R\mathcal{L}_r(M)$, which maps 
links with $\geq r+1$ components to zero, induces the epimorphism
$\mathcal{W}(M)\rightarrow \mathcal{C}_r(M)$. But all relations from immersions with $\geq r+2$ components map trivially because
the links in the corresponding integrability relation have $\geq r+1$ components. 
The relations coming from immersions with $\leq r+1$ components involve
only links with $\leq r$ components. These map trivially into $\mathcal{C}_r(M)$ since the relation submodule $c_r(M)$ of $\mathcal{C}_r(M)$ contains the corresponding Conway triple terms.  
\end{proof}

\begin{remark}
(a)\qua The naturality of the pairings is described by the commuting diagram
$$
\begin{CD}
\mathcal{H}_1(M;R)\otimes \mathcal{C}_{r-1}(M)@>\gamma_{r-1}>> \mathcal{C}_{r-1}(M)\\
@Vid\otimes \iota_{r-1}VV @V\iota_{r-1}VV \\
\mathcal{H}_1(M;R)\otimes \mathcal{C}_r(M)@>\gamma_r>> \mathcal{C}_r(M)
\end{CD}
$$ 

\noindent In particular by theorem 5 the composition
$\gamma_r(x,\ - \ )\circ \iota _{r-1}\circ (z \cdot )$
is trivial on $\mathcal{C}_{r-1}(M)$ for all $x\in \mathcal{H}_1(M)$.
Here $(z \cdot )$ is multiplication by $z$ on $\mathcal{C}_{r-1}(M)$.

\noindent (b)\qua Obviously theorem 1 follows from theorem 5 by taking direct limits. 

\noindent (c)\qua For $r=1$, $\mathcal{C}_1(M)\cong R\hat{\pi}(M)$ and 
$\gamma _1$ is the original Chas-Sullivan pairing.
Moreover by induction $\mathcal{C}_r(M)$ is free if and only
if the pairings $\gamma_1,\gamma_2,\ldots ,\gamma_{r-1}$ are trivial. 

\noindent (d)\qua It has been proved in \cite{CS1} and \cite{K1} that for a hyperbolic $3$-manifold $M$ all pairings $\gamma_r$ (and so $\gamma_{\infty}$) are 
trivial. It is a conjecture (see \cite{K1}) that if $M$ is not hyperbolic then $\gamma_1$ is not trivial and thus $\mathcal{C}(M)$ is not free.  
\end{remark}

\begin{proof}[Proof of theorem 2] Each non-trivial 
$\Theta (\alpha ,\alpha_0,h)\in SR\hat {\pi}(M)$ generates $z$-torsion in the skein module. Thus, if $\mathcal{C}(M)$ is free then all the 
$\Theta $-elements are trivial. But the $\Theta $-elements    
map onto the image of $\gamma _{\infty }$ by a \textit{standard link}
epimorphism:
$$\sigma : SR\hat{\pi}(M)\rightarrow \mathcal{C}(M),$$
which maps $\alpha \in \mathfrak{b}_r(M)$ to some oriented link $K_{\alpha}\in \mathcal{L}_r(M)$ with homotopy classes corresponding to 
$\alpha $ (for details compare \cite{K1}). Thus if $\mathcal{C}(M)$ is free then the image of $\gamma_{\infty }$ is trivial.

Conversely suppose that $\mathcal{C}(M)$ is not free. Choose the 
\textit{smallest} non-negative integer $r$ and a corresponding triple 
$(\alpha ,\alpha_0,h)$ such that $\alpha \in \mathfrak{b}_r(M)$
and $\Theta (\alpha ,\alpha_0,h)\neq 0 \in SR\hat {\pi}(M)$. Let 
$K\in \mathcal{L}_{r-1}(M)$ be the isotopy class of a link with 
homotopy classes of components $\alpha \setminus \alpha_0$.
(Note that by assumption $\mathcal{C}_r(M)$ is free but $\mathcal{C}_{r+1}(M)$ is not free.) It will be shown that $\sigma (\Theta (\alpha ,\alpha_0,h)) 
\notin \mathcal{R}$, so its image in $\mathcal{C}(M)$ is non-trivial.
Note that, since $\Theta (\alpha ,a,h)$ is defined by expanding smoothed links with $r$ components in terms of standard links with respect to $\sigma $, a standard link $K_{\beta }$ with $\beta \in \mathfrak{b}_{r-k}(M)$ appears with a scalar in $\mathbb{Z}[z]$ divisible by $z^k$ in that expansion. 
If $x$ is trivial in $\mathcal{C}(M)$ then it is a linear combination 
of elements $z \Theta (\alpha ',\alpha_0',h')$ with $\alpha '\in \mathfrak{b}_s(M)$ and $s\geq r$. But in such a linear combination 
each $K_{\beta }$ with $\beta \in \mathfrak{b}_{s-k}(M)$ appears 
with a scalar divisible by $z^{s-k+1}$ and $s-k+1>r-k$ for all $s\geq r$. Thus  $\Theta (\alpha ,\alpha_0,h)$ is not contained in 
$\mathcal{R}$.   
\end{proof}

\begin{proof}[Proof of the $\mathbb{Z}$-torsion theorem] 
Let $M=S^2\times S^1$. 
In the following the isomorphisms $\pi_1(M)\cong H_1(M;\mathbb{Z})\cong \mathbb{Z}$ and $S\mathbb{Z}\hat{\pi}(M)\cong S\mathbb{Z}H_1(M;\mathbb{Z})$ are used.

It is easy to show that 
$2zK=0$: Let $K'$ be the $3$-component link formed from $K$ and a meridian curve $S^1\times \{*\}\subset S^2\times S^1$ \textit{linking} both components of $K$. Let $K''$ be the disjoint union of $K$ with an unknot contained in a $3$-ball in $M$. By two applications of the Conway skein relation involving two crossing changes it follows that the difference $K'-K''$ is $2zK$. By using the \textit{belt trick} (handle-slide) it can be seen that the two resulting smoothings are isotopic to $K$. Using the belt trick it also follows that $K'$ and $K''$ are isotopic, which shows that $2zK=0\in \mathcal{C}(M)$.

It will be shown that if  
$nzK=0$ is a relation in $\mathcal{C}(M)$ then $n$ is even.
  
First a new link invariant $lk^2$ in $\mathbb{Z}_2$ is defined for links $L$ in $M$ for which the total homology class $[L]\in H_1(M;\mathbb{Z})$ is even. Let $L_{even}$ (respectively $L_{odd}$) be the sublinks of those components of $L$ with even homology classes (respectively odd homology classes). So $L_{odd}$ has an even number of components. Then $L_{even}$ bounds a (not necessarily oriented) surface. Its mod-$2$-intersection number with $L_{odd}$ is an invariant of links in $M$. In fact, it only has to be checked that this number is not changed by a handle-slide.
But a handle-slide of a component of $L_{even}$ changes this number by the homology class of $L_{odd}$, which is even by assumption. The same argument shows that this number is not changed by a handle-slide on a component of $L_{even}$. Thus $lk^2(L)$ is well-defined. Note that each
crossing change of a component of $L_{even}$ with a component of $L_{odd}$ does change $lk^2(L)\in \mathbb{Z}_2$ while crossing changes among components of $L_{even}$ or $L_{odd}$ do not change $lk^2(L)$. 

In the following assume that $K$ has been chosen to be the standard link for 
the element $1\cdot 1\in \mathfrak{b}_2(M)$.
It follows from the proof of theorem 2 that a relation $nzK=0$ implies that
$$nK=\sum_{h}\lambda_h \Theta (\alpha,\alpha_0,h)\qquad \qquad \qquad (*)$$
where $h$ runs through a finite number of self-homotopies of standard links
$K_{\alpha }$ for $\alpha \in \mathfrak{b}_3(M)$. (By abuse of notation 
the sum is also over a set of those $\alpha$ and corresponding $\alpha_0$.)
This is because relations defined from $\alpha \in \mathfrak{b}_r(M)$ for $r>3$ contain $2$-component links only in higher powers of $z$ and relations from $\mathfrak{b}_2(M)$ contain only elements of $\mathcal{L}_1(M)$.

Now $K$ can appear as a smoothing of a $3$-component link $K_{\alpha }=K'$ only if the homology classes of components of $K'$ are given by $\alpha =1ab\in S\mathbb{Z}H_1(M;\mathbb{Z})$ where $1\in H_1(M;\mathbb{Z})$ and $a+b=1\in H_1(M;\mathbb{Z})$. Without loss of generality we can assume that $a$ is even and $b$ is odd. Let $K_1,K_a,K_b$ denote the corresponding components of $K$. Observe that the total homology class of $K'$ is even and $lk^2(K')$ is defined. Moreover $lk^2(K')$ is changed by crossing changes of $K_a$ with $K_b$ or $K_1$, but remains unchanged by crossing changes of $K_b$ with $K_1$. Now consider a transverse self-homotopy $h$ of $K'$ and the resulting linear combination of smoothings. Only the homology classes of the components of the smoothed links are of importance. In fact, by further applying skein relations only contributions divisible by $z^2$ will appear. Consider only those terms in the linear combination containing terms with homology classes of both components odd. Note that a smoothing of $K_b$ with $K_1$ never can give rise to such a term because it contains an even homology class. If $a\neq 0$ then a smoothing of $K_a$ with $K_1$ gives rise to an (unordered) pair of homology classes $(a+1)b\neq 1\cdot 1\in S\mathbb{Z}H_1(M;\mathbb{Z})$. Assume that along $h$ the smoothings define $n_1(h)$ terms with homology classes $1\cdot 1$ and $n_2(h)$ terms with only odd homology classes but different from $1\cdot 1$, where $n_i(h)\in \{\pm 1\}$ is the sign determined by the crossing change. Since $lk^2(K')$ is well-defined, $n_1(h)+n_2(h)$ is an even number for each self-homotopy $h$. 
 
In the linear combination on the right hand side of $(*)$ all $2$-component links containing even homology classes of components cancel out. Now consider the homomorphism $$R\mathcal{L}(M)\rightarrow \mathbb{Z}\mathcal{L}(M)\rightarrow \mathbb{Z},$$
defined from the coefficient map $R\rightarrow \mathbb{Z}$ and by mapping two-component links to $1\in \mathbb{Z}$, and all other links trivially.  
Consider the image of the right hand linear combination $(*)$ under this homomorphism. 
The resulting finite sum $n'=\sum_{h}\lambda_{h}(n_1(h)+n_2(h))\in \mathbb{Z}$ satisfies $\sum_{h}\lambda_{h}n_2(h)=0$ because of equation $(*)$ and the fact that $K$ has homotopy classes of components $1\cdot 1$. Also $n'$ is even because each 
$n_1(h)+n_2(h)$ is even.
Thus $n=\sum_h\lambda_{h}n_1(h)$ is even.  
\end{proof}

\section{Proof of theorem 3}

It suffices to prove that for $\mathcal{W}(M;\mathbb{Z})$ the pairing $\gamma_u$, which is defined on representative link homotopy classes of links by the Chas-Sullivan procedure according to definition 1 of section 2, is \textit{well-defined}. 
Then $(\mathcal{W}(M;\mathbb{Z}),\alpha_u,\beta_u,\gamma_u)$ is a
$\mathbb{Z}$-deformation of $\gamma $. 
Here $\alpha_u: \mathbb{Z}\mathcal{L}(M)\rightarrow \mathcal{W}(M;\mathbb{Z})$
is the projection, and the natural homomorphism $\mathfrak{h}$ factors through
$\mathcal{W}(M;\mathbb{Z})$ and defines $\beta_u$ (because the homotopy classes of components coincide for $K_{+0}$ and 
$K_{-0}$, respectively $K_{0+}$ and $K_{0-}$). 
The universality follows from the proof below because the integrability relations generate the image in $\mathbb{Z}\mathcal{L}(M)$ of certain differences resulting from applications of the Chas-Sullivan procedure to homologous representative loops in $map(S^1,M)$.
Since these relations have to be true in any geometric $\mathbb{Z}$-deformation of $\gamma $,
an arbitrary deformation factors through the universal one. This means that the homomorphism 
$\alpha_A: \mathbb{Z}\mathcal{L}(M)\rightarrow A$, for a given geometric $\mathbb{Z}$-deformation of $\gamma $, factors through $\alpha_u$. The uniqueness of this factorization follows from the surjectivity of $\alpha_A$.

\begin{remark}
The fact that the universal geometric $R$-deformation of $\gamma $ is defined by the tensor product 
with the universal geometric $\mathbb{Z}$-deformation is immediate from the fact that the integrability relations are integral and from the naturality of all constructions. 
\end{remark}

We first give a precise the definition of $\gamma_u$. 
Let $x\in \mathcal{H}_1(M;\mathbb{Z})$ and $y\in \mathbb{Z}\mathcal{L}(M)$. 

By using the epimorphism $\rho $ defined in the proof of theorem 5 and by linearity, it can be assumed that $x$ is represented by a loop $\ell $ in $map(S^1,M)$ and that $y$ is represented by some embedding $e\in map(\cup_{r-1}S^1,M)$. Then the
\textit{extended} loop $\ell _e: S^1\rightarrow map(\cup_rS^1,M)\subset \widetilde{Map}(S^1,M)$ is defined 
from the union of $\ell $ and $e$ for all parameter values $t$ (compare section 2). 

Assuming weak transversality for $\ell_e$ (compare section 2, lemma 1) an element of $\mathbb{Z}\mathcal{L}(M)$ is defined    
by the signed sum of the link homotopy classes of the smoothings at all singular parameters. The signs are determined by the signs of the double-points (which are defined using the orientation of the loop $\ell $ and thus $\ell_e$, and the orientation of $M$): This is $\pm$ 
if there is a crossing change $K_{\mp} \rightarrow K_{\pm}$ along $\ell_e$ at the 
singular parameter. (Recall that the singular parameters along $\ell_e$ correspond to intersection points of the singular torus in $M$ determined by $\ell $ with the embedding $e$ representing the link homotopy class determined by $e$.)
The image of this linear combination in $\mathcal{W}(M;\mathbb{Z})$ defines
$\gamma_u (x\otimes \bar{y})$, where $\bar{y}$ is the image of $y$ in
$\mathcal{W}(M;\mathbb{Z})$. It has to be shown that this image does not depend on the representatives 
chosen for $x$ and $\bar{y}$. 

First the linear combination is well-defined for a fixed loop $\ell _e$. 
If the embedding $e$ is changed by a link homotopy 
or the loop $\ell $ is changed by a homotopy in the free loop space then there is an induced homotopy of the loop $\ell_e$ in $\widetilde{Map}(S^1,M)$. (A link homotopy of $e$ can delete or produce intersections with the singular torus representing $\ell $.) 
After approximating by a transverse homotopy in Lin's sense (relative to the boundary, compare \cite{K1}, \cite{KL} and \cite{L}), the set of points in $S^1\times I$, for which the image in the mapping space is not in $\mathcal{D}(S^1,M)$, is a properly embedded $1$-complex with vertices of valence $4$ in the interior of $S^1\times I$
and vertices of valence $1$ in the boundary. The \textit{singular parameters} in the boundary components of the annulus (vertices of valence $1$) are joined by properly embedded arcs in the interior, possibly intersecting in points of valence $4$. 

It is easy to see that the two linear combinations according to the boundary components coincide if there are no valence $4$ vertices. In general the homotopy can be decomposed into a finite sequence of \textit{elementary} homotopies containing only a single interior vertex. The  corresponding map $\cup_rS^1\rightarrow M$ is an immersion with exactly two transverse double-points of distinct components. The difference of the linear combinations is (up to a sign) the integrability element determined by that immersion. Because of the very special origin of the homotopy, only immersions in $\mathcal{K}^{2}(M)\cup \mathcal{K}^{3}(M)$ have to be considered. This is because a link homotopy of $e$ (respectively a homotopy of $\ell $) induce a homotopy of
$\ell_e$ for which \textit{all but one} of the components are disjoint during the homotopy. This property can easily be preserved in approximating the homotopy. So the image of the linear combination in $\mathcal{W}(M;\mathbb{Z})$ is not changed under homotopy of the loop $\ell $ or link homotopy of $e$. 

\begin{remark}
For the property of universality of $\gamma _u$ it is important to note that \textit{each} such immersion occurs in a 
homotopy of some loop $\ell_e$ defined as above from a link homotopy or 
homotopy of $\ell $. In fact, only immersions in $\mathcal{K}^{3}(M)$ have to be considered. Then a singular torus together with a null-homotopy and a link homotopy of some embedding $e$ are defined by the neighbourhood of the immersion, which is a family $D^2\rightarrow map(S^1\cup\cup_{r-1}S^1,M)$. 
By a suitable choice of components this family can be induced by a link homotopy of an embedding $e$ or a homotopy of $S^1$-family $\ell $.  
\end{remark}
 
Next consider any loop $\ell $ and and replace $y$ by some integrability element. Then the resulting linear combination defined by the smoothings is an integrability element.
It can be assumed that the two double-points of the immersion defining the element are contained in the complement of possible intersection points of the immersion with the singular torus (defined from $\ell $).
 
We have shown that the \textit{pairing}
$$\bar{\gamma }: \left(\bigoplus_{[f]\in \hat{\pi}(M)}\mathbb{Z}\pi_1(map(S^1,M),f)\right)\bigotimes
\mathcal{W}(M;\mathbb{Z})\rightarrow \mathcal{W}(M;\mathbb{Z})$$
is well-defined.
The pairing is extended over the free abelian group by linear extension. 
Now $\rho $ factors as the composition:
$$\bigoplus_{[f]\in \hat{\pi}(M)}\mathbb{Z}\pi_1(map(S^1,M),f)\rightarrow
\bigoplus_{[f]\in \hat{\pi}(M)}\mathbb{Z}\pi_1(map(S^1,M)/S^1,\bar{f})\rightarrow 
\mathcal{H}_1(M;\mathbb{Z}),$$
using the homomorphism induced by $map(S^1,M)\rightarrow map(S^1,M)/S^1$ (mapping loops to strings), and the Hurewicz homomorphisms. 

It is immediate from definitions that for $x'\in \pi_1(map(S^1,M),f)$ (considered as subset of  $\oplus_{[f]\in \hat{\pi}(M)}\mathbb{Z}\pi_1(map(S^1,M),f)$) the element $\bar{\gamma } (x',y)$ only depends on the image of $x'$ in 
$\pi_1(map(S^1,M)/S^1,\bar{f})$. This is because the link homotopy classes do not
depend on specific choices of basepoints on the loops, from which they result by smoothing, according to the definition of $\bar{\gamma }$.  
Thus the homomorphism $\bar{\gamma }$ factors through the epimorphism $\gamma '$, which is  
defined on $$\bigoplus_{[f]\in \hat{\pi}(M)}\mathbb{Z}\pi_1(map(S^1,M)/S^1,\bar{f}).$$
Let $\rho ': \mathbb{Z}\pi_1(map(S^1,M)/S^1,\bar{f})\rightarrow \mathcal{H}_1(M;\mathbb{Z})$
be the obvious homomorphism factoring $\rho $. Now it is not hard to prove that an element $$\sum_{j}\lambda_jg_j\in \mathbb{Z}\pi_1(map(S^1,M)/S^1,\bar{f})$$ is contained in 
the kernel of $\rho '$ if and only if $\prod_jg_j^{\lambda_j}$ is in the commutator subgroup of 
$\pi_1(map(S^1,M)/S^1,\bar{f})$. But $\gamma '(\prod_jg_j^{\lambda_j},y)= 
\sum_j \lambda_j \gamma '(g_j,y)$ for each $y\in \mathcal{W}(M;\mathbb{Z})$ is immediate 
from the definition of $\bar{\gamma }$. 
Also $\gamma '$ vanishes on commutators by construction.  
So $\bar{ \gamma }(\ \cdot \ ,y)$ factors through $\mathcal{H}_1(M;\mathbb{Z})$ and $\gamma (x,y):=\bar{\gamma}(x',y)$ is well-defined as required. This completes the proof of theorem 3.  \hfill{$\square$} 

\begin{remark} 
(a)\qua It follows from the proof above that if the 
loop $\ell $ is replaced by a loop $S^1\rightarrow map(\cup_rS^1,M)$ for $r\geq 2$ (which replaces the singular torus by a collection of singular tori) then  
we have to divide by integrability relations from $\mathcal{K}^{4}(M)$ too.
The origin of those integrability relations are possible transverse intersections 
among components in the family $\ell $. 

\noindent (b)\qua If we replace in the beginning of our discussion the link homotopy classes by isotopy classes of oriented links in $M$, then also self-intersections of loops in the family $\ell $ have to be taken into consideration. Then immersions with two transverse double points involving precisely one self-crossing would automatically lead to relations
$K_{0+}-K_{0-}=0$ (with the first position indicating the crossing of different components). This is because smoothings of self-crossings are not picked up by the Chas-Sullivan procedure. If applied to immersions for which the double point of distinct components involves an unknotted circle, which intersects the union of
other components in precisely one point, then $K_0=K$ is an immersion, which can be assumed to carry an arbitrary self-crossing (compare the special integrability relations defined in the beginning of section 5).  
So it follows from the integrability relation that arbitrary self-crossings can be changed. This reduces the universality problem of extending $\gamma $ for isotopy of links to link homotopy of links.        

\noindent (c)\qua The argument above could be simplified by considering the epimorphisms
for each $[f] \in \hat{\pi }(M)$:
$$\pi_1(map(S^1,M),f)\rightarrow
\pi_1(map(S^1,M)/S^1,\bar{f})\rightarrow 
H_1(map_f(S^1,M)/S^1;\mathbb{Z}),$$
and noticing that $\mathcal{H}_1(M;\mathbb{Z})\cong \oplus_{[f]\in \hat{\pi}(M)} H_1(map_f(S^1,M)/S^1;\mathbb{Z})$.  
But the proof indicates that $\bigoplus_{[f]\in \hat{\pi}(M)}\mathbb{Z}\pi_1(map(S^1,M),f)$
is the natural object for studying generalized (possibly non-homologically) intersection theory of free loop spaces in order to provide a basis for other interesting deformations. 
\end{remark}

\section{The Chas-Sullivan module and related modules}

Throughout this section let $R=\mathbb{Z}[z]$. We will use integer coefficients for homology and the Chas-Sullivan modules.

Let $\widehat{\mathcal{W}}(M)$ (respectively $\widehat{\mathcal{W}}^+(M)$) be the abelian groups defined by restricting 
integrability relations to those immersions
in $\mathcal{K}^{3}(M)$ (respectively in $\mathcal{K}^{3}(M)\cup \mathcal{K}^4(M)$) for which at least one of the components with a double-point 
is an unknotted circle, which bounds a disk intersecting the other components in precisely one transverse point. These integrability relations are called \textit{special}.    

There are obvious projection homomorphisms:
$$
\begin{CD}
\widehat{\mathcal{W}}(M)@>>> \widehat{\mathcal{W}}^+(M)\\
@VVV @VVV \\
\mathcal{W}(M)@>>> \mathcal{W}^+(M)
\end{CD}
$$

and the same diagram holds for reduced modules. 

For $r\geq 0$ let $\mathcal{L}^r(M):=\mathcal{L}_r(M)\setminus \mathcal{L}_{r-1}(M)$ be the set of link homotopy classes of oriented links in $M$ with $r$ components.
For $\mathcal{Z}(M)$ any of the modules above, 
there is the obvious isomorphism 
$\mathcal{Z}(M)=\bigoplus_{r\geq 0}\mathcal{Z}_r(M)$,
where $\mathcal{Z}_r(M)$ is the quotient of $\mathbb{Z}\mathcal{L}^r(M)$ by the corresponding subgroup of relations. (Note that integrability relations are homogenous with respect to the number of components, and all projections are compatible with this grading.) So there is a commutative diagram as above for each $r\geq 0$.  

By definition $\mathcal{Z}_0(M)\cong \mathbb{Z}$ is generated by the empty link. 
Obviously $\mathcal{Z}_1(M)\cong \mathbb{Z}\hat{\pi}(M)$,
$\mathcal{W}_2(M)\cong \mathcal{W}^+_2(M)$ and 
$\widehat{\mathcal{W}}_2(M)\cong \widehat{\mathcal{W}}^+_2(M)$ holds. 
Also there are induced geometric pairings for $\mathcal{W}^+(M)$ (and for
$\mathcal{W}^+_r(M)$ and $\mathcal{W}_r(M)$ for each $r\geq 0$), which induce structures of 
modules over the Lie algebra $\mathcal{H}_1(M)$ as usual.   

\begin{lemma} 
The abelian group $\widehat{\mathcal{W}}(M)$ 
(respectively  $\widehat{\mathcal{W}}^+(M)$) is isomorphic to the quotient of $\mathbb{Z}\mathcal{L}(M)$ by the submodule generated by skein elements
$$K_+-K_--K_0\sharp H+K_0\sharp U^2,$$
where $U^2$ is the two-component unlink, $H$ is the oriented Hopf link, and $\sharp$ is the oriented connected sum along any of the two components involved in the crossing $K_{\pm}$ (respectively along any component of the link).
\end{lemma}

\begin{proof}[Proof] This is immediate from the definition of the relations. 
\end{proof}

\begin{remark}
\noindent (a)\qua Let $H_-$ be the negatively oriented Hopf link in $M$.   
A special case of the relation in the lemma is:
$$U^2-H_-=H-U^2,$$
so $H_-=2U^2-H$ in $\widehat{\mathcal{W}}(M)$ or $\widehat{\mathcal{W}}^+(M)$.

\noindent (b)\qua The connected sum of a link $K$ with $U^2$ is the disjoint union with the unknot contained in a $3$-ball separated from the link, which in particular does not depend on a choice of component of $K$. 
\end{remark}

\begin{theorem} The inclusion $\mathcal{Z}\mathcal{L}_{\centerdot}(M)\rightarrow R\mathcal{L}_{\centerdot}(M)$ induces the isomorphism $\widehat{\mathcal{W}}_{\centerdot}^+(M)\rightarrow \mathcal{C}_{\centerdot}(M)$.
Thus $\widehat{\mathcal{W}}_{\centerdot}^+(M)\cong \mathcal{W}_{\centerdot}^+(M)\cong \mathcal{C}_{\centerdot}(M)$ as abelian groups.
\end{theorem}

\begin{proof}[Proof] Since the integrability relations (for all possible immersions with two double-points) over $\mathbb{Z}$ map trivially into $\mathcal{C}_{\centerdot}(M)$ the homomorphism
$$\mathcal{Z}\mathcal{L}_{\centerdot}(M)\rightarrow R\mathcal{L}_{\centerdot}(M)\rightarrow \mathcal{C}_{\centerdot}(M)$$ factors through
$\widehat{\mathcal{W}}_{\centerdot}^+(M)$ and defines a homomorphism $\chi : \widehat{\mathcal{W}}_{\centerdot}^+(M)\rightarrow \mathcal{C}_{\centerdot}(M)$.
Now the free abelian group $\mathbb{Z}[z]\mathcal{L}_{\centerdot}(M)$ has the basis given by elements $z^iK$ for $i\geq 0$. First define a homomorphism $\mathfrak{k}$ 
of abelian groups
$\mathbb{Z}[z]\mathcal{L}_{\centerdot}(M)\rightarrow \widehat{\mathcal{W}}_{\centerdot}^+(M)$ 
by induction on $i$. For $i=0$ use the natural 
epimorphism $\alpha : \mathbb{Z}\mathcal{L}_{\centerdot}(M)\rightarrow \widehat{\mathcal{W}}_{\centerdot}^+(M)$. For $i\geq 1$ and $K\in \mathcal{L}_{\centerdot}(M)$ map $z^iK$ to $\alpha \circ (H-U^2)\sharp 
\bar{(\mathfrak{k}}(z^{i-1}K))$, where 
$$(H-U^2)\sharp : \mathbb{Z}\mathcal{L}_{\centerdot}(M)\rightarrow \mathbb{Z}\mathcal{L}_{\centerdot}(M)$$ is defined by the linear extension of a map of link homotopy classes of links $K\mapsto K\sharp H^2-K\sharp U^2$ with connected sum along some chosen component of $K$ and $\bar{\mathfrak{k}}(x)$ is a lift of $\mathfrak{k}(x)$ to $\mathbb{Z}\mathcal{L}_{\centerdot}(M)$ for each $x\in \widehat{\mathcal{W}}_{\centerdot}^+(M)$.
The composition is well-defined and independent of the choice of component for the connected sum. This is because special integrability relations from $\mathcal{K}^{4}(M)$ permit to change between components used in connected sums with Hopf links. Moreover the difference of lifts is a linear combination of integrability relations, which maps to a linear combination of integrability relations under $(H-U^2)\sharp $.  
Since the Conway element $K_+-K_--zK_0$ maps to the special integrability element of the previous lemma 
there is the induced homomorphism $\kappa : \mathcal{C}_{\centerdot}(M)\rightarrow \widehat{\mathcal{W}}_{\centerdot}^+(M)$. 
Now $\kappa \circ \chi =id$ is obvious. The identity $\chi \circ \kappa =id$ follows because an
application of the Conway skein relation to the non-trivial crossing of the \textit{Hopf link crossing} in $H\sharp K-U^2\sharp K$ shows equality with $zK$ in $\mathcal{C}(M)$. Thus claim follows for all $i\geq 1$ by induction. This proves the claim.
\end{proof}

\begin{remark} 
\noindent (a)\qua The reduced Hoste-Przytycki skein module naturally appears
as an additional structure on the module $\mathcal{W}_{\centerdot}^+(M)$. Just define $z:=H-U^2$ and multiplication by $z$ by $zK=K\sharp H-K\sharp U^2$.

\noindent (b)\qua It follows from the $\mathbb{Z}$-torsion theorem that the abelian group $\mathcal{W}_{\centerdot}^+(S^2\times S^1)$ has torsion. Let $K'$ be the link in
$S^2\times S^1$ defined in the proof of the $\mathbb{Z}$-torsion theorem, which is isotopic to $K\cup U$, compare section 3. 
Then, by an application of a special 
skein relation on some immersion in $\mathcal{K}^{3}(M)$, the relation
$2(K\sharp H-K\cup U)=0$ follows in $\mathcal{W}(S^2\times S^1)$. Now $K\sharp H-K\cup U$ maps to
$zK\in \mathcal{C}(M)$ and is non-trivial there by the proof of the $\mathbb{Z}$-torsion
theorem. Thus $K\sharp H-K\cup U$ is a torsion element in the universal abelian group $\mathcal{W}(S^2\times S^1)$ too. Note that, for each oriented $3$-manifold, the image of a link in $\mathcal{W}(M)$ cannot be a torsion element because the composition
$\mathbb{Z}\mathcal{L}(M) \rightarrow \mathcal{W}(M)\rightarrow S\mathbb{Z}\hat{\pi}(M)$ 
is the homomorphism $\mathfrak{h}$ defined in section 3.  
\end{remark}

The study of the difference between $\mathcal{W}(M)$ and $\mathcal{W}^+(M)$ turns out to be subtle and interesting. Note that from the lemma above, $\mathcal{W}(M)$ is generated by standard links corresponding to homotopy classes of components and iterated connected sums of Hopf links with these links. But connected sums with Hopf links in $\mathcal{W}(M)$ \textit{a priori} depend on the choice of components.

The following notion is of importance in the general study of skein modules.

\begin{definition} A generating set $\mathfrak{S}$ of an $R$-module $A$ is called \textit{minimal} if $\mathfrak{S}\setminus \{g\}$ does not generate $A$ for each $g\in \mathfrak{S}$.
\end{definition}

A \textit{minimal} system of generators of $\mathcal{W}(M)$ should take into account that the usual standard links (corresponding to the homotopy classes of components) can be chosen in such a way that duplicate homotopy classes are represented by parallel link components. Then connected sums with Hopf links along different components of the same free homotopy class are \textit{isotopic} links. A connected sum with a Hopf link, for a link with homotopy classes of components $\alpha_1^{n_1}\ldots \alpha_r^{n_r}$ with $\alpha_i\neq \alpha_j$ for $i\neq j$ and $n_i\geq 1$, can be done in $r$ \textit{a priori} different ways. (If $u$ is the conjugacy class of the trivial homotopy class and $u\neq \alpha \in \hat{\pi}(M)$ then the monomial $u^2\alpha \in \mathfrak{b}_3(M)\subset S\mathbb{Z}\hat{\pi}(M)$ represents the homotopy classes of components of a $3$-component standard link with two trivial homotopy classes and one component with homotopy class $\alpha $.) The previous example shows that \textit{a priori} different links are not necessarily different in all $3$-manifolds. A possible choice of a minimal set of generators of $\mathcal{W}(M)$ as abelian group should depend on the $3$-manifold under consideration. 

The next example shows that iterated connected sums contain hidden symmetries with respect to skein relations, which have to be detected. 

\begin{example} Let $1\neq \alpha \in \hat{\pi}(M)$ and let $K$ be a knot with homotopy class $\alpha $. Let $K_{**}'\in \mathcal{K}^{3}(M)$ be defined from $K$ and some unlink of two components, each component intersecting $K$ in precisely one point. Let $K_{**}$ be the 
$4$-component immersion defined by by taking a connected sum of a Hopf link with one of the trivial components of 
$K_{**}$. The integrability relation resulting from this immersion contains a difference
$(K\sharp H)\sharp_1H-(K\sharp H)\sharp_2H$, where $\sharp_1,\sharp_2$ indicate the connected sum along $K$ (respectively the trivial component of $K\sharp H$). The other two terms are of the form
$(K\cup U)\sharp_1 H-(K\cup U)\sharp_2 U$. By iterating this argument it follows that in building up a generating set for $\mathcal{W}(M)$ \textit{chaining} of Hopf links can be avoided, i.e.\ connected sums with $(\ldots ((K_{\alpha }\sharp H)\sharp H) \ldots \sharp H)\sharp H$ can be performed along $K_{\alpha }$ for $\alpha \in \mathfrak{b}(M)$. So for
given $\alpha := \alpha_1^{n_1}\ldots \alpha_r^{n_r}\in \mathfrak{b}_{n_1+\ldots n_r}(M)$ let $r'\leq r$ be the number of $1\leq i\leq r$ with $n_i\geq 2$. Then there are $r(2r'+(r-r'))$ possibly not equivalent $2$-fold connected sums. It seems to be an interesting independent problem to study the involved combinatorics in detail.   
\end{example} 

A possible difference between $\mathcal{W}(M)$ and $\mathcal{W}^+(M)$ can be first detected on the level of $3$-component links in $M$. Because of the observations and the example above it is not surprising that linking numbers play an important role in detecting such a difference.  
In the following see \cite{K5} for all details concerning generalized linking numbers. 

\begin{proof}[Proof of Theorem 4]  For $0\neq a\in H_1(M)$ a $\mathcal{W}$-skein invariant is defined, which is not a $\mathcal{W}^+$-skein invariant. For each link $K$ let
$lk_{a}(K)$ denote the \textit{sum} of the linking numbers in $\mathbb{Q}$
between components with homology classes $a_1,a_2$ for which $a_1+a_2=a$. 
For oriented closed rational homology $3$-spheres the linking numbers in $\mathbb{Q}$ are classical and defined e.g.\ in \cite{K5}.
The linking numbers in $\mathbb{Q}$ are defined by mapping the universal linking numbers (defined in
the commutative ring $\mathcal{R}(M):= \mathbb{Z}[\frac{1}{k}][x_{ij}]$, where $k$ is the torsion order of $H_1(M)$, $1\leq i,j\leq m$, and $m$ is the rank of $H_1(M)$) 
into $\mathbb{Q}$ using a fixed chosen oriented embedding into a closed rational homology $3$-sphere.  
The map $lk_a :\mathcal{L}(M)\rightarrow \mathbb{Q}$ 
extends by linearity to the free abelian group $\mathbb{Z}\mathcal{L}(M)$. 
The only property of $lk_a$ we need here is that it is changed by addition of $\pm 1$ under a crossing
change of components with homology classes summing up to $a$. By definition it does not change under other crossing changes.  

Next let $K_{**}\in \mathcal{K}^{3}(M)$ and let $K_{**}'$ be the immersion $\cup_3S^1\rightarrow M$ containing the double-points.
Let $a_1,a_2,a_3\in H_1(M)$ be three homology classes represented by the components of $K_{**}'$. We can assume without loss of generality that the component appearing at both places $(**)$
has homology class $a_1$. Let $K=K_{**}-K_{**}'$ be the link resulting from $K_{**}$ and let
$K_{\pm 0}'$ (respectively $K_{0 \pm}'$) denote the $2$-component sublinks resulting by resolving the double-points. Then the homology classes of $K_{\pm 0}'$ (respectively $K_{0\pm}'$) are $a_1+a_2,a_3$ (respectively $a_1+a_3,a_2$). The relevant linking numbers between components of $K$, or between a component of $K$ with a component of $K_{0\pm}'$ (respectively $K_{\pm 0}'$), appear in pairs of opposite sign because the crossing change does not change the considered linking numbers and the homology classes.
Now the linking number between the two components of $K_{\pm 0}'$ (respectively $K_{0\pm}'$) is only counted if $a_1+a_2+a_3=a$. But in this case the deficiency in the linking numbers
for $K_{+0}'$ and $K_{-0}'$ cancels the corresponding deficiency for $K_{0+}'$ and $K_{-0}'$.
Thus $lk_a$ defines a homomorphism on $\mathcal{W}^+(M)$ (respectively $\mathcal{W}^+_r(M)$) for each $r\geq 2$.  

Now let $K_{**}\in \mathcal{K}^{4}(M)$ and let $K_{**}'$ be the $4$-component sublink of components containing the double-points. Assume that the homology class of one of the components of $K_{**}'$ is $a\neq 0$ and the other homology classes are trivial. As above the relevant linking numbers counted in $lk_a$ are cancelling in pairs, for the links in the integrability relation, if a component of $K=K_{**}-K_{**}'$ is involved. Note that the three components of each of the links $K_{\pm 0}'$ and $K_{0\pm }'$ have homology classes $a,0,0$. Assume that $a$ is the homology class of one of the two components intersecting in the first place of $(**)$. Then the linking numbers of the smoothed component with any of the two other components in $K_{0\pm}'$ cancel in the integrability relation because of different signs as before. The linking number between the two other components is not counted since $a\neq 0$. 
Finally in the difference $K_{+0}'-K_{-0}'$ the smoothed component has homology class $0$ and the 
linking numbers with the component of homology class $a$ in $K_{+0}'$ (respectively $K_{-0}'$) coincide. So they cancel in the difference. The only remaining linking numbers counted are those between the two components corresponding to the first place in $(**)$. But this precisely creates a difference 
of $1$ in $lk_a$, so $lk_a$ does not vanish on the integrability relation determined by $K_{**}$.            
\end{proof}

\begin{remark} The proof of theorem 4 shows that for $r\geq 3$ the epimorphisms
$$\mathcal{W}_r(M)\rightarrow \mathcal{W}^+_r(M)$$
are not injective under the assumptions $H_1(M)\neq \{0\}$ and $M$ a submanifold of a rational homology $3$-sphere.
\end{remark}

\begin{example}
\noindent (a)\qua Let $K$ be a knot in $M$ with homology class $a \neq 0$.
Then $lk_a((K\sharp H)\cup U)=1$ and $lk_a(K\cup H)=0$.    

\noindent (b)\qua The assumptions of theorem 4 hold e.g.\ for all lens spaces $\neq S^3, S^1\times S^2$ and for all link complements in rational homology $3$-spheres. 
\end{example}

It should be possible to modify the above argument under the weaker assumption that $M$ is not simply connected, by counting suitable linking numbers weighted by elements of $\hat{\pi}(M)$. 
The difficulty is that the operation of smoothing or loop multiplication  
multiplies free homotopy classes \textit{locally}. A transverse double-point of distinct components of
a given immersion defines lifts of the homotopy classes to $\pi_1(M)$ of the two involved components. The free homotopy class of the smoothed component is defined by the 
product at this double-point. It is the \textit{locality} of this multiplication, which prevents that an invariant $lk_{\alpha }$
for $\alpha \in \hat{\pi}(M)$, defined analogous to $lk_a$ above, vanishes in general on
integrability relations of type $\mathcal{K}^{3}(M)$. But a modification of the counting of linking numbers, or an invariant counting the linking numbers 
weighted in $\hat{\pi}^2(M):=\pi_1(M)\times \pi_1(M)/{s.\ c.\ }$, where $s.\ c.\ $ is simultaneous conjugation (compare \cite{C}),   
might detect the precise difference between $\mathcal{W}(M)$ and $\mathcal{W}^+(M)$.

It seems to be an even more difficult problem to detect a possible difference between $\mathcal{W}(M)$ and $\mathcal{W}^+(M)$ for simply connected $M$. The two $4$-component \textit{extended standard links} $(H\sharp H) \cup U$ and $H\cup H$ are the same element in $\mathcal{W}_4^+(S^3)$ and probably are differently in $\mathcal{W}_4(S^3)$. But it is not obvious how to construct a suitable skein invariant detecting the difference. 

\begin{conjecture} \ (a)\qua The projection 
$\mathcal{W}_r(S^3)\rightarrow \mathcal{W}_r^+(S^3)$ is an isomorphism for
$r=3$ but is not an isomorphism for $r\geq 4$. 

\noindent (b) $M$ is simply connected if and only if $\mathcal{W}_3(M)\rightarrow \mathcal{W}_3^+(M)$ is an isomorphism.

\noindent (c)\qua For $M$ a submanifold of a rational homology sphere, a minimal generating set for
$\mathcal{W}(M)$ is defined from all \textit{extended standard links}, i.e.\ all standard links $K_{\alpha }$ for $\alpha \in \mathfrak{b}(M)$ and iterated Hopf links on standard links restricted only by the symmetries described in the example above. It seems to be interesting to consider $\mathcal{W}_{\centerdot}(M)$ as a deformation of $\mathcal{C}_{\centerdot}(M)\cong \mathcal{W}_{\centerdot}^+(M)$ with deformation parameters measuring the different ways of taking connected sums with Hopf links giving rise to \textit{essentially} different links in the $3$-manifold.    

\noindent (d)\qua If $M\neq S^3$ is a Seifert-fibred $3$-manifold then $\mathcal{W}(M)$ has torsion. 
\end{conjecture}

\bibliographystyle{gtart}

\Addresses\recd

\end{document}